\documentclass[review]{elsarticle}

\usepackage{appendix}
\usepackage{amsthm}
\usepackage{lineno}
\usepackage{amsmath,mathtools,esint}
\usepackage{amsfonts}
\usepackage{multirow}% http://ctan.org/pkg/multirow
\usepackage{multicol}
\usepackage{hhline}% http://ctan.org/pkg/hhline
\usepackage[colorlinks,bookmarksopen,bookmarksnumbered,citecolor=red,urlcolor=blue]{hyperref}
\usepackage{verbatim}
\usepackage{caption}
\usepackage{array}
\usepackage{bigints}
\usepackage{soul}
\usepackage{amssymb}
\usepackage{wrapfig}
\usepackage{graphicx}
\usepackage[left=1in, top=1in, right=1in, bottom=1in]{geometry}
\usepackage[utf8]{inputenc}
\usepackage{lipsum}
\usepackage{amsfonts}
\usepackage{graphicx}
\usepackage{epstopdf}
\usepackage{amsopn}
\usepackage{psfrag}
\usepackage[usenames,dvipsnames]{color}
\usepackage[normalem]{ulem}
\usepackage{MACROS/remarks,MACROS/mydef,MACROS/boldfonts,MACROS/mydef_sph,xspace}
\usepackage{pgfplots}
\usepackage{subfiles}
\usepackage[lined,algonl,boxed]{algorithm2e} 
\usepackage{tikz}
\usepackage{subcaption}

\modulolinenumbers[5]

\journal{ArXiv.org}

\begin{document}

\begin{frontmatter}

\title{An efficient algorithm for weakly compressible flows in spherical geometries}

% \address{Radarweg 29, Amsterdam}
% \fntext[myfootnote]{Since 1880.}


\author[mymainaddress]{Roman Frolov}
\ead{frolov@ualberta.ca}

\author[mymainaddress]{Peter Minev\corref{mycorrespondingauthor}}
\ead{pminev@ualberta.ca}

\author[mymainaddress]{Aziz Takhirov}
\ead{takhirov@ualberta.ca}

\address[mymainaddress]{Department of Mathematical and Statistical Sciences, University of Alberta, Edmonton, AB, T6G 2G1, Canada}

\begin{abstract}
{In this paper we present a direction splitting method, combined with a nonlinear iteration, for the compressible Navier-Stokes equations in spherical coordinates.  The method is aimed at solving the equations on the sphere, and can be used 
for a regional geophysical simulations as well as simulations on the entire sphere.  The aim of this work was to develop a method that would work efficiently in the limit of very small to vanishing Mach numbers, and we demonstrate here, using a numerical example, that the method shows good convergence and stability at Mach numbers in the range $[10^{-2}, 10^{-6}]$.  We also demonstrate the effect of some of the parameters 
of the model on the solution, on a common geophysical test case of a rising thermal bubble.  The algorithm is particularly suitable for a massive parallel implementation, and we show below some results demonstrating 
its excellent weak scalability.
}
 \end{abstract}


\end{frontmatter}

%\linenumbers
%%%%%%%%%%%%%%%%%%%%%%%%%%%%%%%%%%%%%%%%%%%%%%%%%%%%%%%%%%%%%%%%%%%%%%%%%%%%%%%%%%%%%%%%%%%%%%%%%%%%%%%%%%%%%%%
%%%%%%%%%%%%%%%%%%%%%%%%%%%%%%%%%%%%%%%%%%%%%%%%%%%%%%%%%%%%%%%%%%%%%%%%%%%%%%%%%%%%%%%%%%%%%%%%%%%%%%%%%%%%%%%
\section{Introduction.}

The main motivation for this study comes from the atmospheric science and oceanography, where reliable dynamical cores for global and local ocean-atmosphere circulations are required to decrease uncertainties in numerical 
weather prediction, ocean circulation, and climate modelling. Despite the rapid advance in numerical methods for atmospheric and oceanic flows, there are still several important challenges remaining in this field. Among them are 
the need to avoid simplifications of the model that may only be valid in certain asymptotic limits, improve the efficiency and increase the accuracy and resolution of the computations while maintaining stability, consistently couple 
ocean and atmosphere models together, and many others. 

Very generally speaking, the atmospheric and ocean models can be divided into two large classes: hydrostatic and non-hydrostatic. In the first case the flow in the vertical direction is "homogenized" based on  
the assumption that the spherical shell domain of the flow is very thin as compared to its size in the other directions, while in the second case the model explicitly involves in the system of equations the balance of mass, momentum, 
and energy in the vertical direction as well. Of course, the non-hydrostatic models are  the most comprehensive models of the atmosphere and the oceans, and naturally, they are based on the 3D compressible Navier-Stokes 
equations. This system can be further enhanced by introducing equations accounting for the moisture content and pollutants in the atmosphere, salinity of the oceans, etc., thus enhancing the accuracy of the modelling effort, but 
also increasing the complexity of the PDE system, and the required computational effort for their numerical approximation.  Therefore, until the appearance of the modern parallel computer systems, the hydrostatic approach was prevailing. However, in the last two decades, the modelling community clearly more and more often employs the complete 3D models.  But even the 3D computational cores are often based on some simplified models 
 like the compressible/incompressible  Euler equations or the incompressible Boussinesq equations.  The computational reason for considering the Euler equations is usually the fact that they can be efficiently discretized by means 
 of fully explicit schemes while the incompressible models avoid the numerical difficulties related to the treatment of the weakly compressible flows.  Excellent summaries of the present models in meteorology and climatology can be
 found in  \cite{MKMMKVGHJ2016} and in oceanology - in \cite{ocean_circ_2019}. These articles contain also numerous references to various modelling and computational efforts in these areas, and therefore we refer the reader 
 to these references for obtaining a more complete picture of the state-of-the-art techniques for resolution of such models. Our effort in the present study is focussed on the development of an efficient parallel algorithm for the solution of the 3D weakly compressible 
 Navier-Stokes equations in a spherical shell, using semi- or fully implicit schemes in time and a finite difference approximation in space.  The main guidelines for the design of the algorithm were: (i) minimization of the computational 
 effort; (ii) maximization of its parallel performance.  Although we do not claim that we completely achieved those goals, we do believe that this study is a step in the right direction that will eventually lead to the ability to resolve the 
 coupled ocean-atmosphere system with satisfactory resolution on the scale of the entire Earth.  The discretization methodology in this study is applied to the compressible Navier-Stokes equations with the so-called Stiffened Gas 
 (SG) equation of state in spherical coordinates. The advection terms are written in a non-conservative form, { very often used in case of incompressible flows, and therefore } more appropriate, in our opinion, for weakly compressible flows.  This set is further discretized implicitly using the 
 Linearized-Block-Implicit (LBI) direction splitting scheme based on the second order Douglas splitting \cite{Douglas1962}, that allows for a stable and cheap integration in time (see \cite{frolov2019efficient} and the references therein).  In addition, it also 
 allows for a very efficient parallelization if this time discretization is combined with a finite difference approximation in space.  We use a staggered Marker-and-Cell (MAC) grid instead of the centred non-staggered discretization
 that is commonly employed in case of higher Mach numbers.  {  Since this study is focussed on low Mach number flows, this choice  prevents,  without the introduction of any stabilization terms,  the appearance of node-to-node pressure oscillations, typical for incompressible simulations on colocated grids. As it has been already demonstrated in case of incompressible flow, the combination of a direction splitting time discretization with a MAC finite difference discretization in space, yields an algorithm with good stability properties and excellent parallel performance (see for example \cite{doi:10.1137/140975231} and \cite{minevadi}).}

The rest of the paper is organized as follows. We describe the details of the proposed algorithm in Section \ref{3Algorithm}. Numerical experiments are described in Section \ref{3num}, and we provide some concluding remarks in Section \ref{3conclusion}.

\clearpage
\section{Formulation and its discretization}
\label{3Algorithm}
\subsection{A non-conservative formulation in spherical coordinates} 
The flows in the atmosphere and the ocean occur at extremely low to moderate values of the Mach number, and therefore no shock waves are observed in the solution. { Furthermore, it is desirable that a method 
for low Mach number compressible flows also resolves reasonably well the limit of incompressible flows by setting the Mach number to a very low value. Since, arguably the most popular formulation in the incompressible case is the formulation in which the advection terms are written in a non-conservative form, we first re-write the conservative set of compressible equations, widely used for atmospheric modelling in such a form.  The system is formulated in terms of the primitive variables pressure $p$, velocity $\bu$, and temperature $T$, again widely used in the incompressible regime. In case of dry, stratified air, it reads (see  \ref{AppendixA} for a detailed derivation of this system based on the conservative equation set for atmospheric modelling  provided for example in \cite{MKMMKVGHJ2016} \footnote{The system contains a dimensionless parameter, the Prandtl number $Pr$, as well as dimensional parameters.  This is a somewhat unusual setting, introduced in \cite{MKMMKVGHJ2016} but for the sake of consistency with this publication we keep it here.}). : }
\begin{align}
%\begin{equation}
%\begin{split}
 \label{nonconservmass}
 \frac{\partial T}{\partial t} + \bu \cdot \nabla T + (\gamma-1)T \nabla \cdot {\bu} - &\frac{(\gamma-1)T}{p+\pi_{\infty}}\nabla \cdot \left( \frac{\mu c_p}{Pr} \nabla T \right) \\&- \frac{(\gamma-1)T}{p+\pi_{\infty}}\nabla {\bu} : \hat{\bsigma} = 0,\notag   \\
%\end{split}
%\end{equation}
%\begin{equation}\label{Cnonconservmom}
%\begin{split}
\label{nonconservmom}\frac{\partial {\bu}}{\partial t} + {\bu} \cdot \nabla {\bu} + \frac{1}{\rho} \nabla p - \frac{1}{\rho} \nabla \cdot \hat{\bsigma} +& {\bg} + 2 ({\bu} \times {\bomega}) =  0,  \\
%\end{split}
%\end{equation}
%\begin{equation}\label{Cnonconservener}
%\begin{split}
\label{nonconservener}\frac{\partial p}{\partial t} + {\bu} \cdot \nabla p  +\gamma (p+\pi_{\infty}) \nabla \cdot {\bu}    - &(\gamma-1)\nabla \cdot \left( \frac{\mu c_p}{Pr} \nabla T \right) \\&- (\gamma-1)\nabla {\bu} : \hat{\bsigma} = \notag 0.
%\end{split}
%\end{equation}
\end{align}
where ${\bomega}$ is the rotational velocity of the Earth, $\hat{\bsigma}$ is the viscous stress tensor given by $$\hat{\bsigma} = \mu \left[ \left(\nabla {\bu} + (\nabla {\bu})^T \right) - \frac{2}{3} (\nabla \cdot {\bu}) \hat{I} \right],$$ ${\bg}$ is the sum of the true gravity and the centrifugal force, $\displaystyle c_p$, $\displaystyle c_v$, $\mu$, Pr, $\displaystyle \gamma = \frac{c_p}{c_V}$, $\pi_{\infty}$ are constant for each material, and the density $\rho$ is given by the Stiffened Gas Equation of State (\cite{sgeos}):
\begin{equation}\label{nonconserveos}
\begin{split}
\rho = \frac{p + \pi_{\infty}}{c_V (\gamma-1)T}.
\end{split}
\end{equation}

\begin{figure}[htbp]
\centering
\begin{minipage}[b]{0.4\linewidth}
\centering
\includegraphics[width=\textwidth]{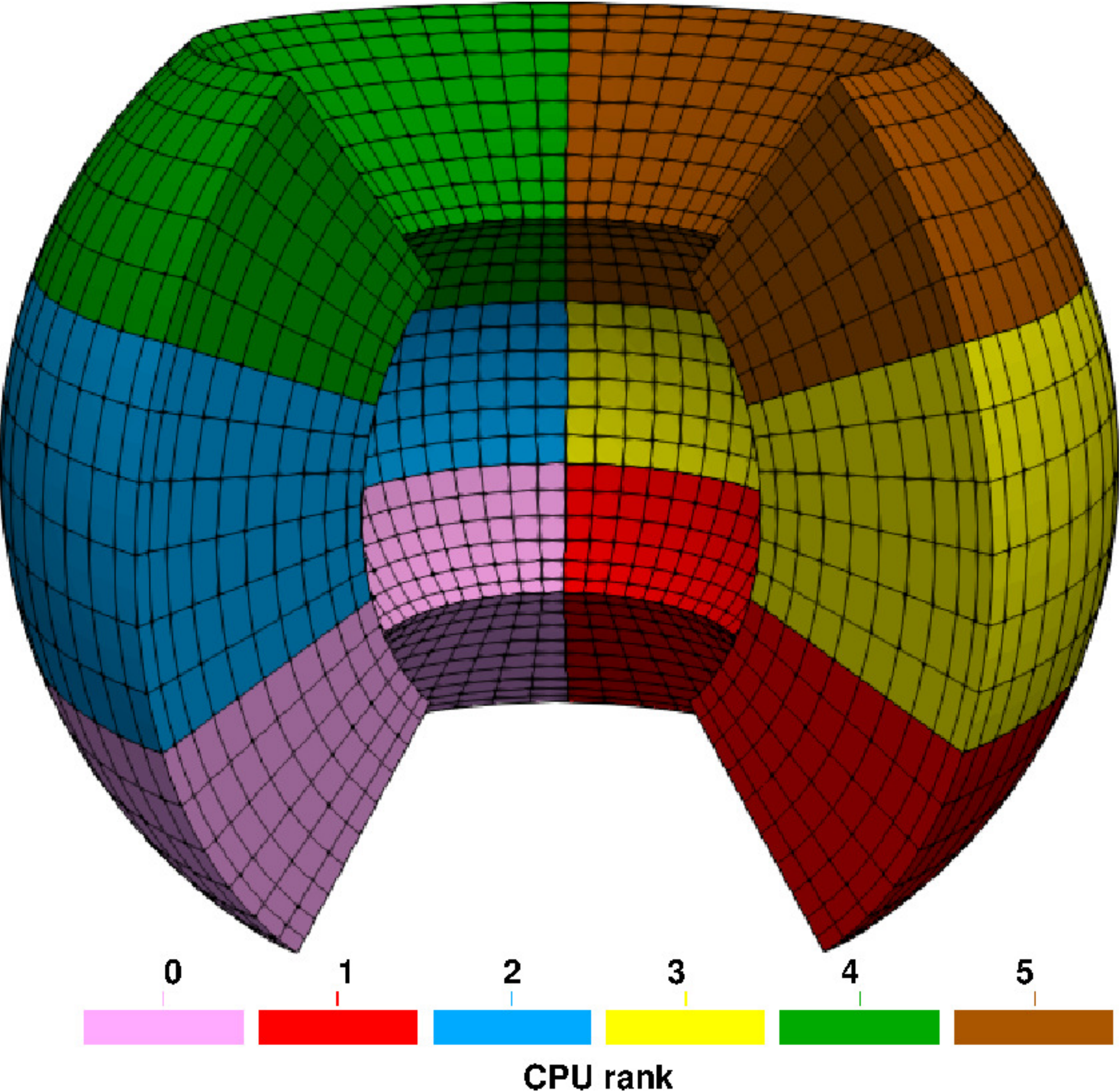}
%\subcaption{Yin and Yang grids for CPU distribution $1\times3\times2$}
\end{minipage}
\hspace{0.1cm}
\begin{minipage}[b]{0.4\linewidth}
\centering
\includegraphics[width=\textwidth]{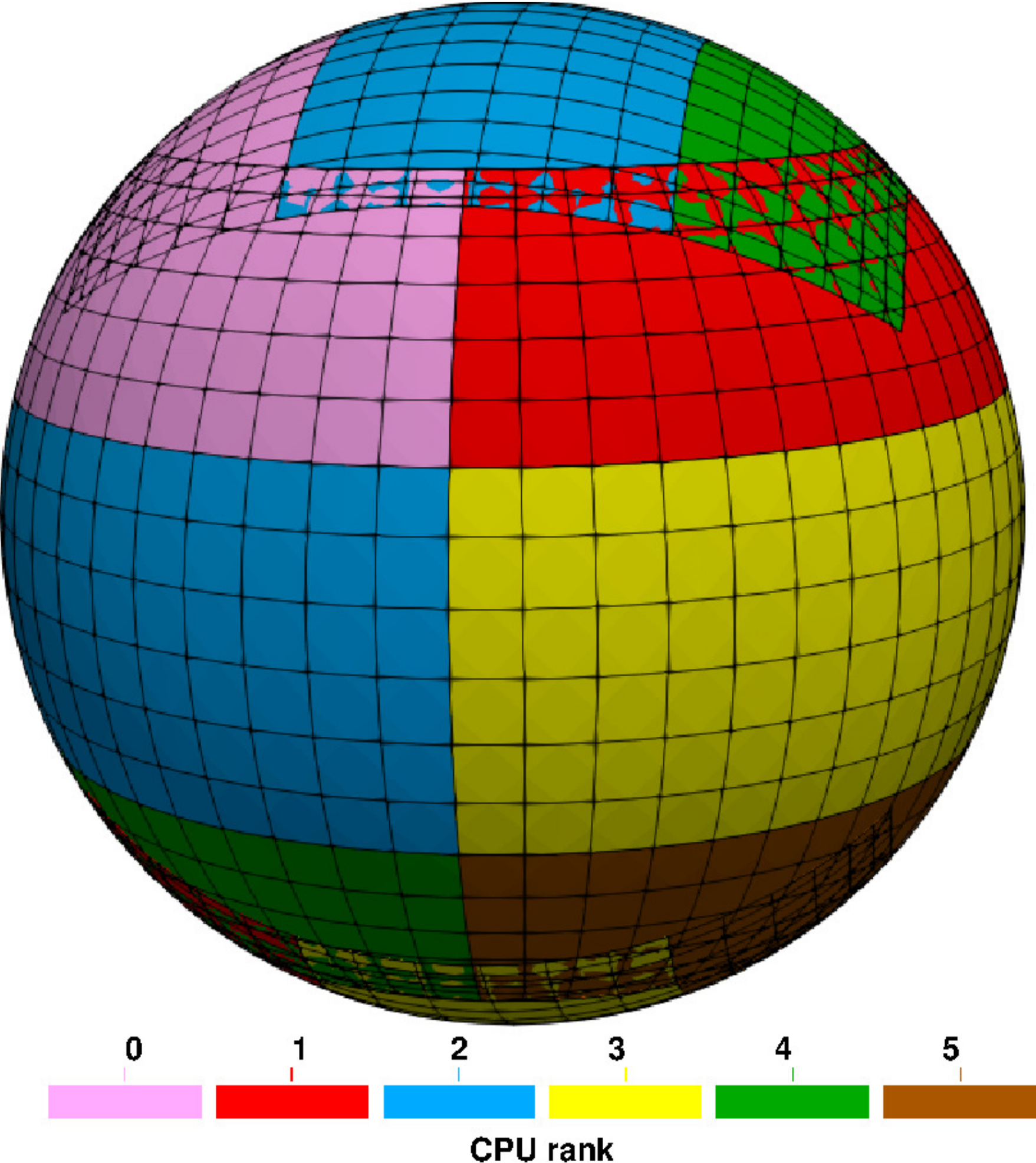}
%\subcaption{Yin-Yang grid with overlap}
\end{minipage}
\caption{Yang (left) and Yin-Yang (right) grids.  Each subgrid is further decomposed into blocks for a parallel implementation corresponding to a CPU distribution $1\times3\times2$.}
\label{fig:YYgrid}
\end{figure}

{ The ultimate goal of this study is to develop a technique that can be applied to the modelling of geophysical flows on a scale of several hundred kilometres to the scale of the entire Earth's atmosphere (and the oceans, after a proper coupling is developed).  Therefore, the speed and parallel efficiency of such a technique is of a paramount importance to us.  With respect to parallel efficiency, perhaps the best choice for a spatial discretization
is to use structured grids, properly adapted around topographical obstacles (see \cite{AKM12} for details).  Since the domain is a-part-of- or a full spherical shell, such a grid choice predetermines the use of a spherical transform.
Perhaps the major problem of such a choice is that the use of a uniform grid size in spherical coordinates would mean a very non-uniform grid size in physical coordinates, that even vanishes around the poles.  To overcome this problem
the spherical transform can be combined with an overlapping domain decomposition approach  for the grid on the entire sphere.  One possible choice for a domain decomposition is provided by the so-called Yin-Yang grids (see \cite{KT2004} and \cite{tfm20}). Such a choice allows for a very easy Cartesian decomposition of the available processors, and therefore makes the parallel implementation quite straightforward. A graphical representation of such
a domain decomposition together with a possible Cartesian processor distribution is given in figure \ref{fig:YYgrid}.  Since, as it will become clear in the next section, the discrete operators are not positive, the design
of the domain decomposition iteration requires serious efforts.  Therefore, in this study we make the first step and consider the system  (\ref{nonconservmass})-(\ref{nonconservener}) in  a part of a spherical shell given by:
$$ \Omega: = \left\lbrace \left(r,\theta,\phi \right) \in [r_1, r_2] \times \left[\theta_1, \theta_2 \right] \times 
\left[ \phi_1, \phi_2 \right]  \right\rbrace. $$}

In order to simplify the notations, we denote the vector of unknowns as $\bU = [p,u_r,u_{\theta},u_{\phi},T]^T$, the gravity vector as $\bG = [0,\bg^T,0]^T$, and combine all the components of the differential operators in the corresponding directions into the $\pmb{D_{r}}(\bU)$, $\pmb{D_{\theta}}(\bU)$, and $\pmb{D_{\phi}}(\bU)$ operators, and all the mixed derivatives, derivates in staggered directions, and other terms not suitable for implicit treatment by the direction splitting approach into the $\pmb{D_{M}}(\bU)$ operator (see (\ref{dr}),(\ref{dtheta}), (\ref{dphi}), and (\ref{dm}) for definitions of these operators). Then the system (\ref{nonconservmass})-(\ref{nonconservener}) can be written in a compact form as (see  \ref{AppendixB} for details):
\begin{equation}\label{chapter3main}
\frac{\partial \bU}{\partial t} + \pmb{D_{r}}(\bU) \bU + \pmb{D_{\theta}}(\bU) \bU + \pmb{D_{\phi}}(\bU) \bU + \pmb{D_{M}}(\bU)\bU + \bG = 0.
\end{equation}

\subsection{Discretization}
The time discretization starts with the Crank-Nicolson time discretization, combined with a Picard nonlinear iteration (see e.g. \cite{picard}, Chapter 3), that yields the following the semi-discrete version of (\ref{chapter3main}):
\begin{equation}\label{chapter3mainsd}
\begin{split}
\frac{\bU^{n+1,k+1} - \bU^n}{\tau} +& \frac{1}{2}\pmb{D}^{n+\frac{1}{2},k} \bU^{n+1,k+1}+\frac{1}{2}\pmb{D}^{n+\frac{1}{2},k}\bU^n +\\& \frac{1}{2}\pmb{D_M}^{n+\frac{1}{2},k} \bU^{n+1,k}+\frac{1}{2}\pmb{D_M}^{n+\frac{1}{2},k}\bU^n+ \bG = 0,
\end{split}
\end{equation}
where $$\pmb{D}(U) = \pmb{D_{r}}(\bU)+\pmb{D_{\theta}}(\bU)+\pmb{D_{\phi}}(\bU),$$ $$\pmb{D}^{n+\frac{1}{2},k} = \pmb{D}\left(\frac{\bU^{n+1,k}+\bU^n}{2}\right),$$ $\tau$ is the time-step, $n$ refers to the time level, and $k=0, 1, \dots, K$ refers to the iteration level. As usual, we set $ \bU^{n+1,0}= \bU^{n}$ and $ \bU^{n+1}= \bU^{n+1,K}$. For the sake of brevity we will skip the superscripts of operators in the remainder of the paper, assuming that $\pmb{D}=\pmb{D}^{n+\frac{1}{2},k}$.
{ Note that the spherical transform introduces some mixed derivatives as well as other terms that, if treated implicitly, would make the coupling between the equations of the semi-discrete system
very strong and hard to handle at the linear algebra level.  Therefore, we discretize these terms, all multiplied by the operator $ \pmb{D_{M}}$,   explicitly with respect to the iteration level.}

In order to reduce further the computational effort, we approximate (\ref{chapter3mainsd}) by a Douglas-type (see \cite{Douglas1962}) direction-wise factorization  that can be written as:
\begin{equation}\label{splitting1}
\begin{split}
\left(I +\frac{\tau}{2}\pmb{D_{r}} \right)\left(I +\frac{\tau}{2}\pmb{D_{\theta}} \right)\left(I +\frac{\tau}{2}\pmb{D_{\phi}} \right) \left( \bU^{n+1,k+1} - \bU^n\right) =\\
-\tau \pmb{D} \bU^n - \tau \bG - \frac{\tau}{2}\pmb{D_{M}}\bU^{n+1,k} - \frac{\tau}{2}\pmb{D_{M}}\bU^n.
\end{split}
\end{equation}

{Since the operators $\pmb{D_{r}}, \pmb{D_{\theta}},$ and $\pmb{D_{\phi}}$ are not positive, it is of a paramount importance for the stability and accuracy of such a direction splitting approximation
to reduce the splitting error given by:
\begin{equation*}
\begin{split}
& \pmb{ER}\left( \bU^{n+1,k+1} - \bU^n\right) = \frac{\tau^2}{4} \left(\pmb{D_{r}} \pmb{D_{\theta}} + \pmb{D_{r}} \pmb{D_{\phi}} + \pmb{D_{\theta}} \pmb{D_{\phi}}+ \frac{\tau}{2} \pmb{D_{r}} \pmb{D_{\theta}} \pmb{D_{\phi}} \right) \left( \bU^{n+1,k+1} - \bU^n\right) =\\
& \left(I +\frac{\tau}{2}\pmb{D_{r}} \right)\left(I +\frac{\tau}{2}\pmb{D_{\theta}} \right) \left(I +\frac{\tau}{2}\pmb{D_{\phi}} \right) \left( \bU^{n+1,k+1} - \bU^n\right) - \left(I + \frac{\tau}{2}\pmb{D} \right)\left( \bU^{n+1,k+1} - \bU^n\right)
\end{split}
\end{equation*}
To our knowledge, the only stable approach for such a splitting error reduction is to add the error at the previous iteration level to the right-hand-side of (\ref{splitting1}).  Note, that any attempt
to use a higher order extrapolation of the error, involving more previous iteration levels, leads to a destabilization of the iteration.
}

After some rearrangement, we obtain the following factorized direction splitting iteration to be solved until convergence:

\begin{equation}\label{splitting2}
\begin{split}
&\left(I +\frac{\tau}{2}\pmb{D_{r}} \right)\left(I +\frac{\tau}{2}\pmb{D_{\theta}} \right)\left(I +\frac{\tau}{2}\pmb{D_{\phi}} \right) \left( \bU^{n+1, k+1} - \bU^{n+1,k}\right) =\\
& -\left(I + \frac{\tau}{2}\pmb{D} \right)\left( \bU^{n+1, k }- \bU^n\right)
-\tau \pmb{D} \bU^n - \tau \bG - \frac{\tau}{2}\pmb{D_{M}}\bU^{n+1,k} - \frac{\tau}{2}\pmb{D_{M}}\bU^n.
\end{split}
\end{equation}
Note that if the number of iterations $K=1$, the resulting time discretization is equivalent to the Douglas scheme, that formally has a $\mathcal{O}(\tau^2)$ splitting error.  Furthermore, the splitting error reduction approach described above has a computational complexity of the same order as the nonlinear Picard iterations without the splitting error reduction.

The system (\ref{splitting2}) can be solved as a sequence of three one-dimensional problems:
\begin{align} 
&\left( I + \frac{\tau}{2}\pmb{D_r} \right) \left( \bxi^{n+1} - \bU^{n+1,k} \right) = -\left(I + \frac{\tau}{2}\pmb{D} \right)\left( \bU^{n+1,k} - \bU^n\right) \label{spl1} \\ 
&\notag -\tau \pmb{D} \bU^n - \tau \bG - \frac{\tau}{2}\pmb{D_{M}}\bU^{n+1,k} - \frac{\tau}{2}\pmb{D_{M}}\bU^n,\\
&\left( I + \frac{\tau}{2}\pmb{D_{\theta}} \right) \left( \beeta^{n+1} - \bU^{n+1,k} \right) = \label{spl2}\bxi^{n+1} -  \bU^{n+1,k},\\
&\left( I + \frac{\tau}{2}\pmb{D_{\phi}} \right) \left( \bU^{n+1,k+1} - \bU^{n+1,k} \right) = \label{spl3}\beeta^{n+1} -  \bU^{n+1,k}
\end{align}
where $\bxi^{n+1}$, $\beeta^{n+1}$, and $\bU^{n+1, k+1}$ are subsequent approximations of the exact solution at $t^{n+1}$.
{ Since the discrete operators involved in the splitting are not symmetric and positive definite, the analysis of this scheme is hard and beyond the goals of this paper.  If  the iteration (\ref{spl1})-(\ref{spl3}) converges at each time step, the resulting solution would be the same as the solution of a Crank-Nicolson 
discretization of the original system (\ref{chapter3main}).  However, proving convergence of the iteration, or stability of the scheme at $K=1$ is far from trivial.  Our numerical experience shows that if one or only a few 
iterations are performed than the overall scheme is stable under a CFL-like condition. Since the diffusivity coefficients for atmospheric flows are negligibly small, one might think that it would be better to simply discretize 
the whole system fully explicitly.  However, explicit methods for compressible flows are subject to the CFL stability condition $\displaystyle \tau \leq \frac{d}{|\lambda_{max}|}$, where $\tau$ is the time-step, $d$ is  the minimum spatial grid size, and $\lambda_{max}$ is the fastest characteristic wave speed, which can be written in terms of flow speed $u$ and sound velocity $c$ as $\lambda_{max} = u \pm c$ (see \cite{cordier}). In dimensionless form this condition becomes:
\begin{equation}
\label{nondimcfl}
\bar{\tau} \leq M_0 \frac{d}{\max|M_0 \bar{u} \pm \bar{c}|},
\end{equation}
where the dimensionless quantities are marked with bars (see \cite{cordier} for details). It is clear from (\ref{nondimcfl}) that the time-step restriction becomes more and more severe as the Mach number $M$ decreases, leading to larger computational time (over-resolution in time).  In the limit of a zero Mach number such a scheme would be unconditionally unstable.
The artificial compressibility method for incompressible flow overcomes this by penalizing the incompressibility constraint, but even then  the fully explicit discretization of the Boussinesq system (with an explicit treatment
of the pressure) is stable only under a condition that the time step is of order of $d^2$ (see \cite{Glow_2003}, chapter IV).
Since the purpose of the present method is to eventually be 
applicable to oceanic flows too, we give preference to the present iterative approach that makes it more flexible, with the computational expense being a constant multiple of the one of a fully explicit scheme.
}

\begin{figure}[ht]
\centerline{
\scalebox{0.2}{\includegraphics[]{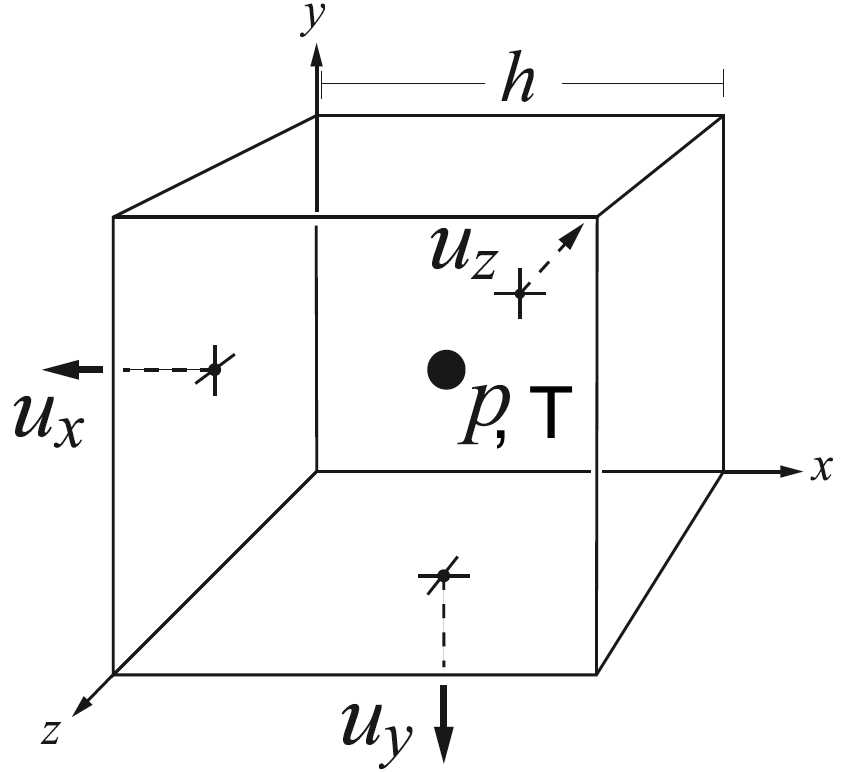}}
}
\caption{MAC stencil with indication of the positions of the variables $p, \bu, T$.}
\label{fig:mac}
\end{figure}

{ Similarly to our previous work on incompressible flow (see \cite{doi:10.1137/140975231}, \cite{minevadi}) we apply a staggered MAC finite difference discretization that yields an {\em inf-sup} stable discretization at zero Mach number.  { A Cartesian MAC cell, including the position of the different variables is shown in figure \ref{fig:mac}. }Since the incompressible flow is a limiting case of the weakly compressible regime we expect that this property of the MAC discretization would prevent the appearance of artificial pressure oscillations at low Mach numbers as well.   Note that the first derivatives are approximated with standard central difference on the MAC stencil, and no artificial stabilization is added.
This means that all the stability for hyperbolically dominated cases is provided by the level of implicitness, that is controlled by the convergence of the iterations described above.  Additionally, the standard finite difference discretization on Cartesian grids is ideally suited to be combined with direction splitting methods since it yields a block-tridiagonal system in each spatial direction after the discretization of the problems (\ref{spl1})-(\ref{spl3}). The solution of these systems can be performed by the block-tridiagonal extension of the Thomas algorithm (see \cite{fletcher1}, Volume 1, pp.188-189), and the operation count of the resulting method is only slightly higher than a fully explicit approach. The parallel implementation of the Thomas algorithm using the Schur complement technique with a non-overlapping domain decomposition, as described in \cite{minevadi}, can be easily extended for the block-tridiagonal version of the linear solver. Weak scalability results for the method can be found in \cite{frolov2019efficient}.}

\section{Numerical tests.}
\label{3num}
The numerical experiments presented below confirm the accuracy of the proposed scheme in a wide range of Mach numbers ($M \in [10^{-6},10^{-1}]$), the correct behaviour of the numerical solution in the incompressible limit, and excellent parallel performance of the method.

\subsection{Well-prepared manufactured solution.}
The following manufactured solution has been used to verify the implementation and study the properties of the algorithm:
\begin{equation}  \label{man_sol}
\begin{aligned}
\rho =& \rho_0 = 1, \\
p = & p_0 + u_0^2 \left(1 + \sin (5t) + \cos^2(\pi r)\cos^2 (4 \phi) \cos^2 (4\theta) \right),\\
u_r =& \frac{u_0 (1+ \sin t)}{2 r^2} + \frac{u_0^2}{c_0} \left(1 + \sin (4t) + \sin\left( r^2 \right) \cos^3(\theta)\sin^2(\phi) \right), \\
u_{\theta} =& \frac{u_0 (1+ \cos (3t+2))}{2 \sin \theta} + \frac{u_0^2}{c_0} \left(1 + \sin (t) + \cos^3\left( r^2 \right) \cos^2(\theta)\sin^3(\phi) \right),\\
u_{\phi} =& \frac{u_0 (1+ \sin (6+t))}{2 } + \frac{u_0^2}{c_0} \left(1 + \cos (2+t) + \cos\left( r \right) \sin^3(\theta)\sin^2(\phi) \right),\\
T =& \frac{p}{c_v (\gamma-1)}.
\end{aligned}
\end{equation}
Note that this solution provides well-prepared initial data, i.e. it has the correct scaling with respect to the Mach number. Thus, it can be used to study the behaviour of the scheme in the incompressible ($M_0 \rightarrow 0$) limit. Indeed, since the characteristic density $\rho_0=1$, the characteristic sound speed becomes $\displaystyle c_0 = \sqrt{\gamma p_0} \sim \sqrt{p_0}$, and the characteristic Mach number is equal to $\displaystyle M_0 = \frac{u_0}{c_0} \sim \frac{u_0}{\sqrt{p_0}}$. Then, non-dimensionalized pressure is given by $$\displaystyle \tilde{p} = \frac{p}{p_0} = \tilde{p}_0 + M_0^2 \tilde{p}_2(r,\theta,\phi,t),$$ and the non-dimensionalized divergence $$\frac{\nabla \cdot {\bu}}{u_0} \sim M_0,$$ which is in agreement with the results from \cite{guil1}.

The governing equations are modified by the inclusion of source terms, computed using the manufactured solution. In all the tests presented here $p_0 = 6250$, $\displaystyle \gamma = \frac{c_p}{c_v} = 1.6$, $\mu=1$, $Pr=1$, ${\bomega}={\b0}$, $\bg=\mathbf{0}$, and the domain is given by: $$ \Omega: = \left\lbrace [1, 2] \times \left[ \frac{\pi}{4}, \frac{3\pi}{4}\right] \times \left[ \frac{\pi}{4}, \frac{7\pi}{4} \right]  \right\rbrace. $$
{ Two iterations are performed at each time step, initialized by $\bU^{n+1,0} = \bU^n$.}
 Dirichlet boundary conditions are imposed for the velocity components at all the boundaries, and zero Neumann conditions are used for pressure and temperature (satisfied exactly by the manufactured solution). Different values of $u_0$ may be chosen to study the properties of the algorithm at different characteristic Mach numbers.

First, we examine the space convergence properties at different values of $M_0$. Figure \ref{fig:conv} { presents the $l^2$ norm (also known as the discrete $L^2$ norm) of the error in the  pressure (figures in the left column) and the $\phi$-velocity component (figures in the right column) vs. the grid size $d$, {at $\tau = 10^{-5}$}, for Mach numbers $M_0=10^{-2}$, $M_0=10^{-4}$, and $M_0=10^{-6}$, respectively. Overall, the error in the $l^2$  norm is of a second order.  Although it seems that it slightly deviates from second order at $M_0=10^{-6}$ and small grid sizes, this is due to the round off errors since the error reaches very small values close to $e^{-30}$.}  Although not shown here, the velocity components in $r$ and $\theta$ directions, and temperature exhibit similar convergence rates.

%%%%%%%%%%%%%%%%%%%%%%%%%%%
\pgfplotstableread{ch_3_data/spacem2.dat}{\smt}
\pgfplotstableread{ch_3_data/spacem4.dat}{\smf}
\pgfplotstableread{ch_3_data/spacem6.dat}{\sms}
%%%%%%%%%%%%%%%%%%%%%%%%%%%
\begin{figure}[htp] 
	\begin{subfigure}[b]{.5\linewidth}
		\centering
		\begin{tikzpicture}[scale=0.7]
		\begin{axis}[legend pos=north west,
		xlabel={$\ln{(d)}$},
		ylabel={$\ln{(l^2\text{error})}$}
		]
		\addplot [black,mark = diamond] table [x={lndiam}, y={lnpl2}] {\smt};
		\addlegendentry{$\ln{(L^2\text{error}(p))}$}
		\addplot [blue,dashed] table [x={lndiam}, y={lndiammod}] {\smt};
		\addlegendentry{2nd-order slope}
		\end{axis}
		\path
		% ([shift={(-5\pgflinewidth,-5\pgflinewidth)}]current bounding box.south west)
		([shift={( 5\pgflinewidth, 5\pgflinewidth)}]current bounding box.north east);
		\end{tikzpicture}	
		\subcaption{}
		\label{fig:convpa}
	\end{subfigure}
	\begin{subfigure}[b]{.5\linewidth}
		\centering
		\begin{tikzpicture}[scale=0.7]
		\begin{axis}[legend pos=north west,
		xlabel={$\ln{(d)}$},
		ylabel={$\ln{(l^2\text{error})}$}
		]
		\addplot [black,mark = diamond] table [x={lndiam}, y={lnwl2}] {\smt};
		\addlegendentry{$\ln{(L^2\text{error}(u_{\phi}))}$}
		\addplot [blue,dashed] table [x={lndiam}, y={lndiammodw}] {\smt};
		\addlegendentry{2nd-order slope}
		\end{axis}
		\path
		% ([shift={(-5\pgflinewidth,-5\pgflinewidth)}]current bounding box.south west)
		([shift={( 5\pgflinewidth, 5\pgflinewidth)}]current bounding box.north east);
		\end{tikzpicture}	
		\subcaption{}
		\label{fig:convwa}
	\end{subfigure}
	
	\begin{subfigure}[b]{.5\linewidth}
		\centering
		\begin{tikzpicture}[scale=0.7]
		\begin{axis}[legend pos=north west,
		xlabel={$\ln{(d)}$},
		ylabel={$\ln{(l^2\text{error})}$}
		]
		\addplot [black,mark = diamond] table [x={lndiam}, y={lnpl2}] {\smf};
		\addlegendentry{$\ln{(L^2\text{error}(p))}$}
		\addplot [blue,dashed] table [x={lndiam}, y={lndiammod}] {\smf};
		\addlegendentry{2nd-order slope}
		\end{axis}
		\path
		% ([shift={(-5\pgflinewidth,-5\pgflinewidth)}]current bounding box.south west)
		([shift={( 5\pgflinewidth, 5\pgflinewidth)}]current bounding box.north east);
		\end{tikzpicture}	
		\subcaption{}
		\label{fig:convpb}
	\end{subfigure}
		\begin{subfigure}[b]{.5\linewidth}
		\centering
		\begin{tikzpicture}[scale=0.7]
		\begin{axis}[legend pos=north west,
		xlabel={$\ln{(d)}$},
		ylabel={$\ln{(l^2\text{error})}$}
		]
		\addplot [black,mark = diamond] table [x={lndiam}, y={lnwl2}] {\smf};
		\addlegendentry{$\ln{(L^2\text{error}(u_{\phi}))}$}
		\addplot [blue,dashed] table [x={lndiam}, y={lndiammodw}] {\smf};
		\addlegendentry{2nd-order slope}
		\end{axis}
		\path
		% ([shift={(-5\pgflinewidth,-5\pgflinewidth)}]current bounding box.south west)
		([shift={( 5\pgflinewidth, 5\pgflinewidth)}]current bounding box.north east);
		\end{tikzpicture}	
		\subcaption{}
		\label{fig:convwb}
	\end{subfigure}

\begin{subfigure}[b]{.5\linewidth}
	\centering
	\begin{tikzpicture}[scale=0.7]
	\begin{axis}[legend pos=north west,
	xlabel={$\ln{(d)}$},
	ylabel={$\ln{(l^2\text{error})}$}
	]
	\addplot [black,mark = diamond] table [x={lndiam}, y={lnpl2}] {\sms};
	\addlegendentry{$\ln{(L^2\text{error}(p))}$}
	\addplot [blue,dashed] table [x={lndiam}, y={lndiammod}] {\sms};
	\addlegendentry{2nd-order slope}
	\end{axis}
	\path
	% ([shift={(-5\pgflinewidth,-5\pgflinewidth)}]current bounding box.south west)
	([shift={( 5\pgflinewidth, 5\pgflinewidth)}]current bounding box.north east);
	\end{tikzpicture}	
	\subcaption{}
	\label{fig:convpc}
\end{subfigure}
\begin{subfigure}[b]{.5\linewidth}
	\centering
	\begin{tikzpicture}[scale=0.7]
	\begin{axis}[legend pos=north west,
	xlabel={$\ln{(d)}$},
	ylabel={$\ln{(l^2\text{error})}$}
	]
	\addplot [black,mark = diamond] table [x={lndiam}, y={lnwl2}] {\sms};
	\addlegendentry{$\ln{(L^2\text{error}(u_{\phi}))}$}
	\addplot [blue,dashed] table [x={lndiam}, y={lndiammodw}] {\sms};
	\addlegendentry{2nd-order slope}
	\end{axis}
	\path
	% ([shift={(-5\pgflinewidth,-5\pgflinewidth)}]current bounding box.south west)
	([shift={( 5\pgflinewidth, 5\pgflinewidth)}]current bounding box.north east);
	\end{tikzpicture}	
	\subcaption{}
	\label{fig:convwc}
\end{subfigure}
	\caption{  Manufactured solution (\ref {man_sol}), $log$-$log$ plots of the $l^{2}$ norm of the pressure (left column) and $u_{\phi}$ (right column)  errors  vs. the grid size $d$;  $t=10^{-3}$, $\tau = 10^{-5}$,  Mach numbers: $M=10^{-2}$ (top row), $M=10^{-4}$ (middle row), and $M=10^{-6}$ (bottom row).}
	\label{fig:conv}
\end{figure}
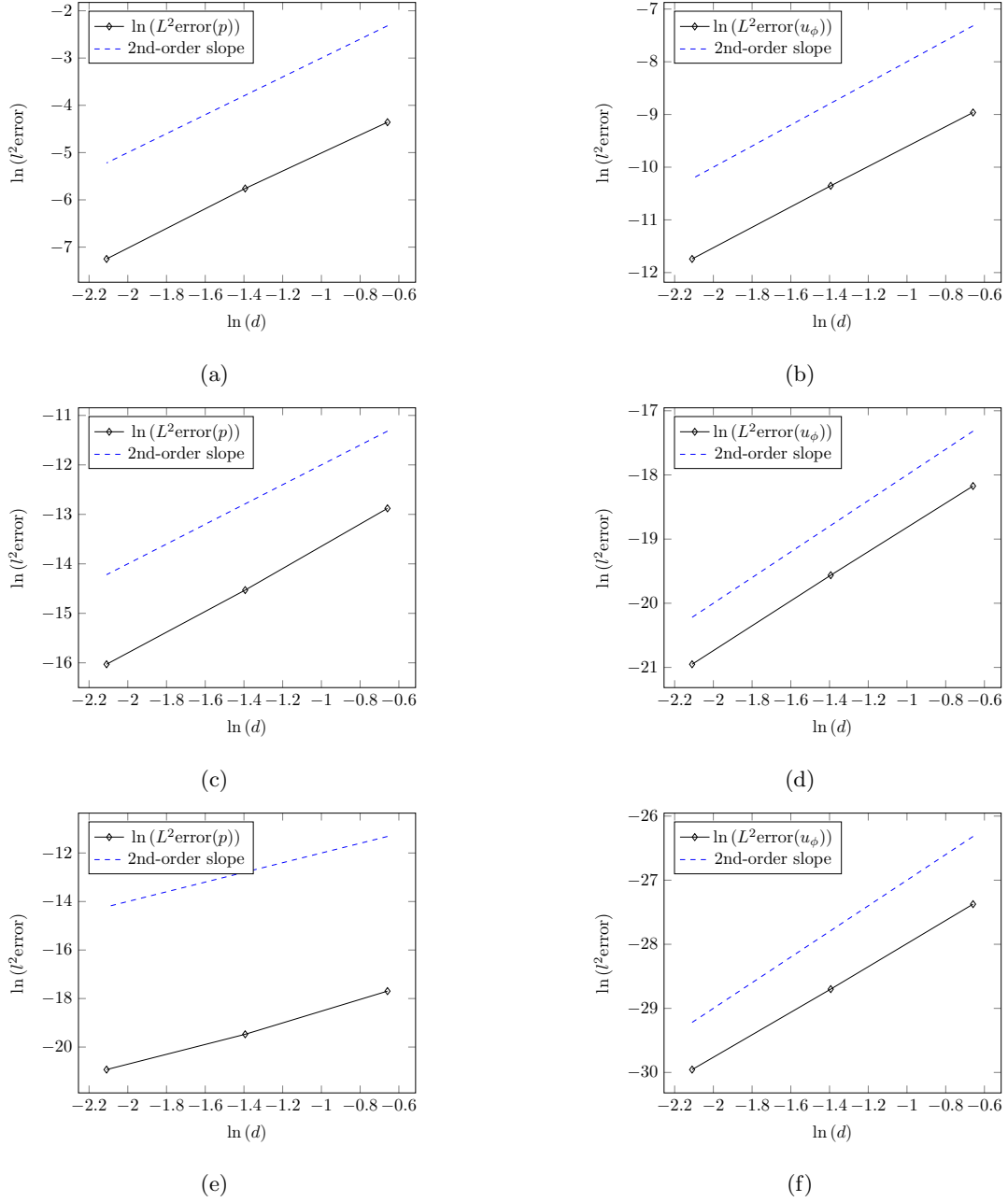

\begin{table}
	\centering
	\begin{tabular}{ |p{1.5cm}||p{2.5cm}|p{2.5cm}|p{2.5cm}|  }
		\hline
		$\tau$ & $M_0=10^{-2}$ & $M_0=10^{-4}$ & $M_0=10^{-6}$\\
		\hline
		$2 \cdot 10^{-3}$ & ($1.8$, $1.7$) &($3.2$, $3.2$)&($3.2$, $3.2$)\\ 
		$1 \cdot 10^{-3}$ &  ($1.9$, $1.9$)& ($2.7$, $2.6$)& ($2.7$, $2.5$)\\
		$5 \cdot 10^{-4}$  & ($1.9$, $1.9$) & ($1.6$, $1.5$)& ($1.4$, $1.5$)\\
		\hline
	\end{tabular}
	\caption{Order of time convergence for $(p,T)$ computed using the inverse Richardson extrapolation approach.}
	\label{tab:timept}
\end{table}

\begin{table}
	\centering
	\begin{tabular}{ |p{1.5cm}||p{2.5cm}|p{2.5cm}|p{2.5cm}|  }
		\hline
		$\tau$ & $M_0=10^{-2}$ & $M_0=10^{-4}$ & $M_0=10^{-6}$\\
		\hline
		$2 \times 10^{-3}$ & ($2.3$, $2.9$, $2.6$) &($3.2$, $2.5$, $2.9$)&($2.5$, $2.9$, $2.4$)\\ 
		$1 \times 10^{-3}$ &  ($2.4$, $2.7$, $3.2$)& ($2.4$, $2.8$, $2.6$)&  ($2.4$, $2.8$, $2.3$) \\
		$5 \times 10^{-4}$  & ($2.1$, $2.4$, $2.3$) & ($2.2$, $2.6$, $2.1$)& ($2.2$, $2.6$, $2.1$)\\
		\hline
	\end{tabular}
	\caption{Order of time convergence for $(u_r,u_{\theta},u_{\phi})$ computed using the inverse Richardson extrapolation approach.}
	\label{tab:timeu}
\end{table}

Next, we follow \cite{tcr} to estimate the order of temporal accuracy using the following time convergence rate (TCR) estimate:
\begin{equation*}
TCR(u_i,\tau) = \log_2 \left[ \frac{||u_i^{\tau} - u_i^{\frac{\tau}{2}}||}{||u_i^{\frac{\tau}{2}} - u_i^{\frac{\tau}{4}}||} \right]
\end{equation*}
Due to the form of the $TCR$, spatial discretization errors cancel (i.e., the leading order truncation error is $\displaystyle const \cdot \tau^l$, where $l$ is the order of accuracy in time), and $TCR \approx l$. This approach is a form of "inverse Richardson extrapolation", and allows for  temporal benchmarking of the algorithm without extreme grid refinement in $3$D. The $TCR$ parameters for pressure and temperature corresponding to different values of $\tau$ and $M_0$ are listed in Table \ref{tab:timept}, while Table \ref{tab:timeu} gives the $TCR$ values for the velocity components. The results demonstrate that the temporal convergence is consistent with the theoretically expected second order.

Although the stability of the scheme has not been studied rigorously, based on our numerical experience the scheme is conditionally stable.
 { One should bear in mind that if the iteration is fully converged, the solution would be identical to a solution of the Crank-Nicolson discretization of the entire system that, at least in the linear case, is stable if the 
eigenvalues.of the discrete operator have non-positive real parts.   However, since in most practical situations the dissipative terms are negligibly small, convergence may require many iterations or the iteration may not be convergent if the time step is too large.  Our numerical experience is that in order to preserve the stability of the scheme, it is advisable that the Courant number is not much larger than one.
It also suggests that the stability restriction does not depend much on the Mach number, in the range $10^{-6} \leq M_0 \leq 10^{-1}$. 
Some dependence on the scaling of the problem, in particular on the value of $p_0$, was observed. For the same value of the Mach number, if $p_0$ is very large, $p_0 > 10^4$, the 
algorithm suffers of very large errors, since the errors in computing pressure derivatives are scaled with this constant.  This can destabilize the algorithm and therefore, the pressure must be rescaled {\em a priori}.}

It is well known, that in the limit to the incompressible regime, some numerical schemes can exhibit noticeable and artificial acoustic waves. Table \ref{tab:inc} provides the maximum norm of the relative pressure fluctuations ($\Delta p = \frac{p-p_0}{p_0}$) after $n$ time steps { with $\tau = 10^{-3}$}, for different values of $M_0$. Theoretically predicted order of magnitude ($\Delta p \sim \mathcal{O}(M_0^2)$) is well preserved by the scheme. The method does not introduce 
noticeable artificial acoustic waves ($\mathcal{O}(M_0)$ pressure fluctuations) even for extremely low values of $M_0$.  
\begin{table}
	\centering
	\begin{tabular}{ |p{1cm}||p{2cm}|p{2cm}|p{2cm}|p{2cm}|p{2cm}|  }
		\hline
		$n$  & $M_0=10^{-2}$ & $M_0=10^{-3}$ & $M_0=10^{-4}$& $M_0=10^{-5}$& $M_0=10^{-6}$\\
		\hline
		1 & $3.2 \cdot 10^{-4}$ &$3.2 \cdot 10^{-6}$&$3.2 \cdot 10^{-8}$&$3.2 \cdot 10^{-10}$&$3.9 \cdot 10^{-12}$\\ 
		50 &  $3.6 \cdot 10^{-4}$& $3.6 \cdot 10^{-6}$&  $3.7 \cdot 10^{-8}$&$4.2 \cdot 10^{-10}$&$9.6 \cdot 10^{-12}$ \\
		100  & $4.0 \cdot 10^{-4}$ & $4.0 \cdot 10^{-6}$& $4.1 \cdot 10^{-8}$&$4.9 \cdot 10^{-10}$&$1.4 \cdot 10^{-11}$\\
		\hline
	\end{tabular}
	\caption{Maximum norm of relative pressure variations ($\Delta p = \frac{p-p_0}{p_0}$) after $n$ time steps for different values of $M_0$. $\tau = 10^{-3}$, grid diameter: $0.12$.}
	\label{tab:inc}
\end{table}

{ \subsection{Rising thermal bubble}
In this section we present results for a common benchmark problem for a dry atmosphere, initially in a hydrostatic equilibrium, the temperature being suddenly perturbed inside a bubble of radius $R$ at $t=0$.
Such benchmarks are usually specified in terms of the potential pressure and temperature:
$$\pi = \left(\frac{p}{p_{00}}\right)^{R/c_p}, \quad \Theta = \frac{T}{\pi},$$
$p_{00}$ being the pressure at the bottom of the domain, usually taken to be $p_{00}=10^5 Pa$.
The initial pressure distribution $\pi_0$ is taken to be at equilibrium with gravity, i.e., $\displaystyle \frac{d \pi_0}{dx} = -\frac{g}{c_p \Theta_0}$,
the initial potential temperature being a constant, $\Theta_0=300 K$.  We consider here two previously proposed benchmark problems.  

The first thermal bubble benchmark is due to \cite{bf2002}, and we denote it by {\em Thermal 1}. The domain is a rectangle with a  horizontal size of $20 km$, 
and a height of $10 km$.  The initial temperature perturbation is given by:
\begin{equation}
 \Delta \Theta = \left\{  \!\! \ \begin{array}{ll} \displaystyle 2 \cos^2 \left(\frac{\pi L}{2} \right),   & \mbox{ if } L \leq R ,\\
 0, & \mbox{ otherwise }. \end{array} \right.
 \label{temp_pert}
 \end{equation}
with $\displaystyle L=\sqrt{\left(\frac{x-x_c}{R}\right)^2 +\left(\frac{y-y_c}{R}\right)^2}$, 
$R=2 km$ being the initial bubble radius, and $x_c=2 km, y_c=10 km$.
We adapted this case to a full 3D setting since our code is designed for a 3D spherical shell geometry.  The domain that we used is a piece of a spherical shell with $r_1=6371 km$ (the approximate Earth radius), $r_2=6381 km$,
$\theta\in [\pi/2-10/6371, \pi/2 + 10/6371], \phi \in [\pi -10/6371, \pi + 10/6371]$.  The other parameters characterize dry air at about $300K$: $Pr=0.71$,  $\displaystyle c_p=1000 \frac{J}{kg K}, c_v = 713  \frac{J}{kg K},
R= 287  \frac{J}{kg K}, \pi_{\infty}=0 Pa, {\bomega}={\b0}s^{-1}, \bg=9.80665 \times  (1,0,0)^T \frac{m}{s^2}$.  The viscosity, unless otherwise specified, is equal to the viscosity of dry air $\displaystyle \mu = 1.846 \times 10^{-5} \frac{kg}{m s}.$  The initial perturbation is the same as in (\ref{temp_pert}), however, $\displaystyle L=\sqrt{\left(\frac{x-x_c}{R}\right)^2 +\left(\frac{y-y_c}{R}\right)^2+\left(\frac{z-z_c}{R}\right)^2}$, with $x_c=2km, y_c=z_c=0 km$.
Because the Earth's radius is huge compared to the domain size, for all practical purposes the domain is a parallelepiped.  Because the bubble in our case is axisymmetric, the resistance of the external fluid is smaller and therefore
the bubble accelerates faster and the velocity becomes large and requires the use of very small time steps or time step adaptivity.  This is why we needed to stop the simulation after $500 s$, earlier than the final time of $1000 s$ of the simulation in  \cite{bf2002}.  The shape of the bubble also significantly differs, the edges not being rolled up as in the case of the 2D simulation in  \cite{bf2002}, Fig. 1.   Since this testcase is posed in a domain of
a more realistic size, we use it to verify the influence of the various parameters on the solution.  The thermal isolines at $t=500 s$ on a grid of $100 \times 200 \times 200$ MAC cells, for various values of the time step $\tau$, dynamic viscosity $\mu$ and number of iterations $K$ are presented in figure \ref{fig:bub1}.  Clearly the viscosity has almost no effect on the results in the range of $1.846 \times 10^{-5} kg/m.s$ (viscosity of air, see figure \ref{fig:bub_b})  to $1.846 \times 10^{-2} kg/m.s$ (figure \ref{fig:bub_c}), and needs to be increased $10^6$ times in order to have some noticeable smoothing effect on the solution (see figure \ref{fig:bub_d}).  This conclusion  justifies to some extend
the use of constant artificial viscosity that is very common in atmospheric simulations.  Of course, this effect is demonstrated on one example only and therefore is not a proof for the validity of such practices.
The decrease of the time step below $1 s$ (figure \ref{fig:bub_b}) to $\tau=0.2 s$ (figure \ref{fig:bub_a}) has no effect on the solution. 
However, a significant increase of the time step to $\tau=50$ using only two iterations per time step (figure \ref{fig:bub_e}) has a dramatic effect on the accuracy of the method and qualitatively  changes  the solution,
dissipating it very rapidly, and noticeably decreasing the speed of the bubble.  The increase of the iterations to fifty (figure \ref{fig:bub_f}) compensates somewhat these effects, however, the solution is still qualitatively very 
different from the solution at $\tau=1$.  The conclusion that can be drawn is that the choice of the time step and number of iterations is crucial for the quality of the simulation, and this suggests the use of time adaptive
algorithms.

The second thermal bubble benchmark is due to \cite{gkc2013}, and we denote it by {\em Thermal 2}.  It is fully 3D, however, the geometry is still Cartesian, the domain being much smaller, of size of $1 km$, and the bubble radius being $R=250 m$. Our domain in spherical coordinates is $[6.371 \times 10^6, 6.37126] \times 10^6] \times [\pi/2-0.5/6371, \pi/2 + 0.5/6371] \times  [\pi -0.5/6371, \pi + 0.5/6371]$, and all the other parameters being the same as
in the previous example with $\mu = 1.846 \times 10^{-5} kg/m.s$ (viscosity of air). The perturbation in the potential temperature is given by \footnote{A personal communication with Dr. Giraldo revealed a typo in their definition of the temperature perturbation in page B1176.  The correct value of $\theta_c$ is 0.25.}:

\begin{equation*}
 \Delta \Theta = \left\{  \!\! \begin{array}{ll} \displaystyle \frac{1}{4} \left( 1 + \cos \left(\pi L \right) \right),   & \mbox{ if } L \leq R ,\\
 0, & \mbox{ otherwise }. \end{array} \right.
 \end{equation*}
 where $\displaystyle L=\sqrt{\left(\frac{x-x_c}{R}\right)^2 +\left(\frac{y-y_c}{R}\right)^2+\left(\frac{z-z_c}{R}\right)^2}$, with $ x_c=260 m, y_c=z_c=0 km$.
 
 We resolved the problem on a grid of $200 \times 200 \times 200$ MAC cells, using time steps $0.25 s$ and $1 s$, changing the number of iterations in the second case from 2 to 10.  The isolines of the potential temperature 
 perturbation (the actual potential temperature minus $300 K$) are presented in figure \ref{fig:bub2}.  The time step has again a dramatic effect on the solution, the larger time step making it more smoothed.  This effect must be due to the
 splitting error since the increase of the number of iterations eliminates it to a large extent.  The results in the first row of the figure are qualitatively comparable to the result in figure 4.1 of \cite{gkc2013}.  More detailed 
 comparison is not reasonable since the equations solved in this paper are different (compare \ref{nonconservmass} to equations 2.1a-2.1c of \cite{gkc2013}). Furthermore, the present algorithm does not use any artificial 
 stabilization unlike the method in \cite{gkc2013} (note that the dynamic viscosity and thermal diffusivity used here are the actual physical values for dry air).

\begin{figure}[htp] 
	\begin{subfigure}[b]{.5\linewidth}
		\centering
			\includegraphics[width=\textwidth]{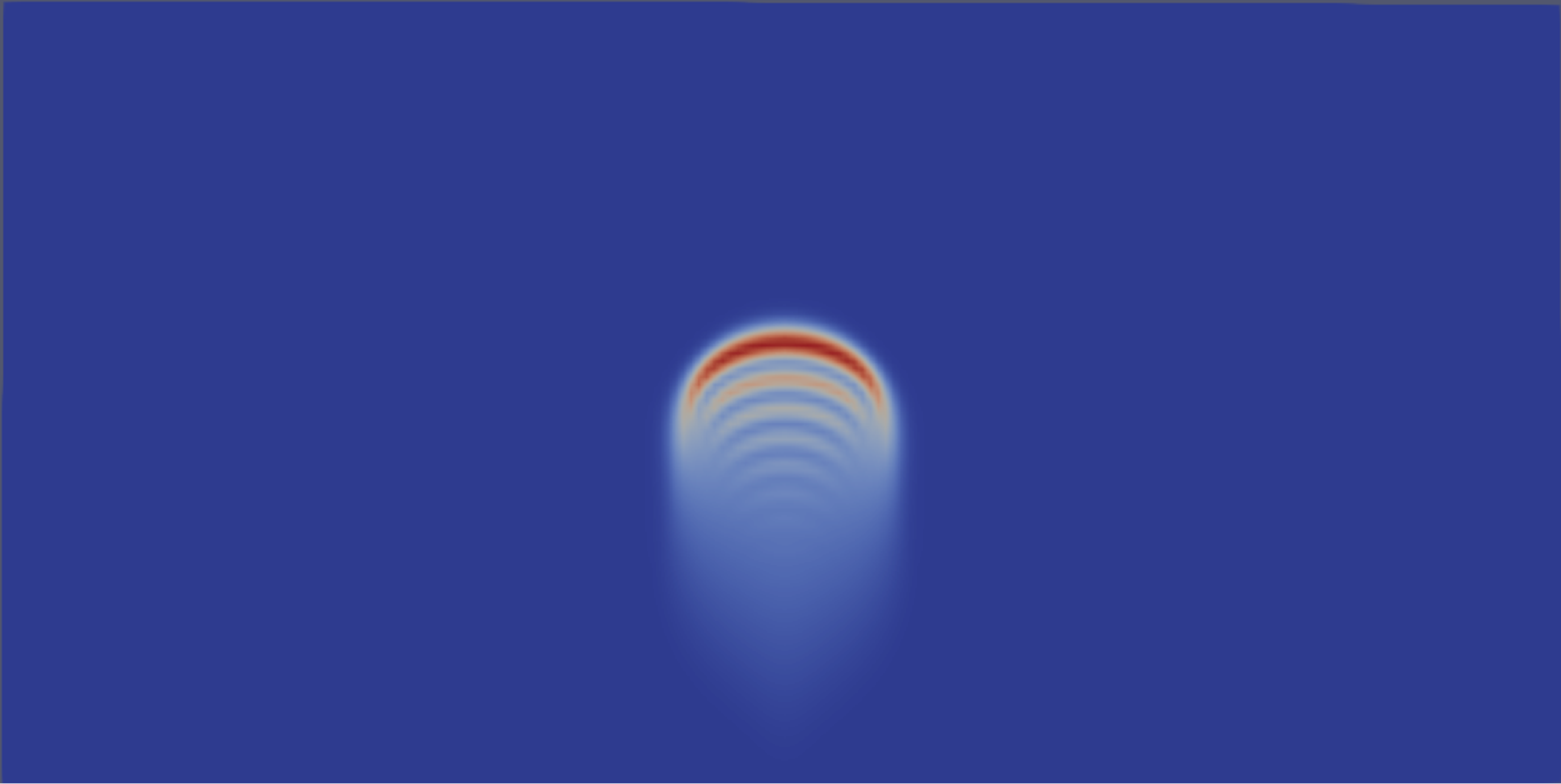}		
		\subcaption{}
		\label{fig:bub_a}
	\end{subfigure}
	\begin{subfigure}[b]{.5\linewidth}
		\centering
			\includegraphics[width=\textwidth]{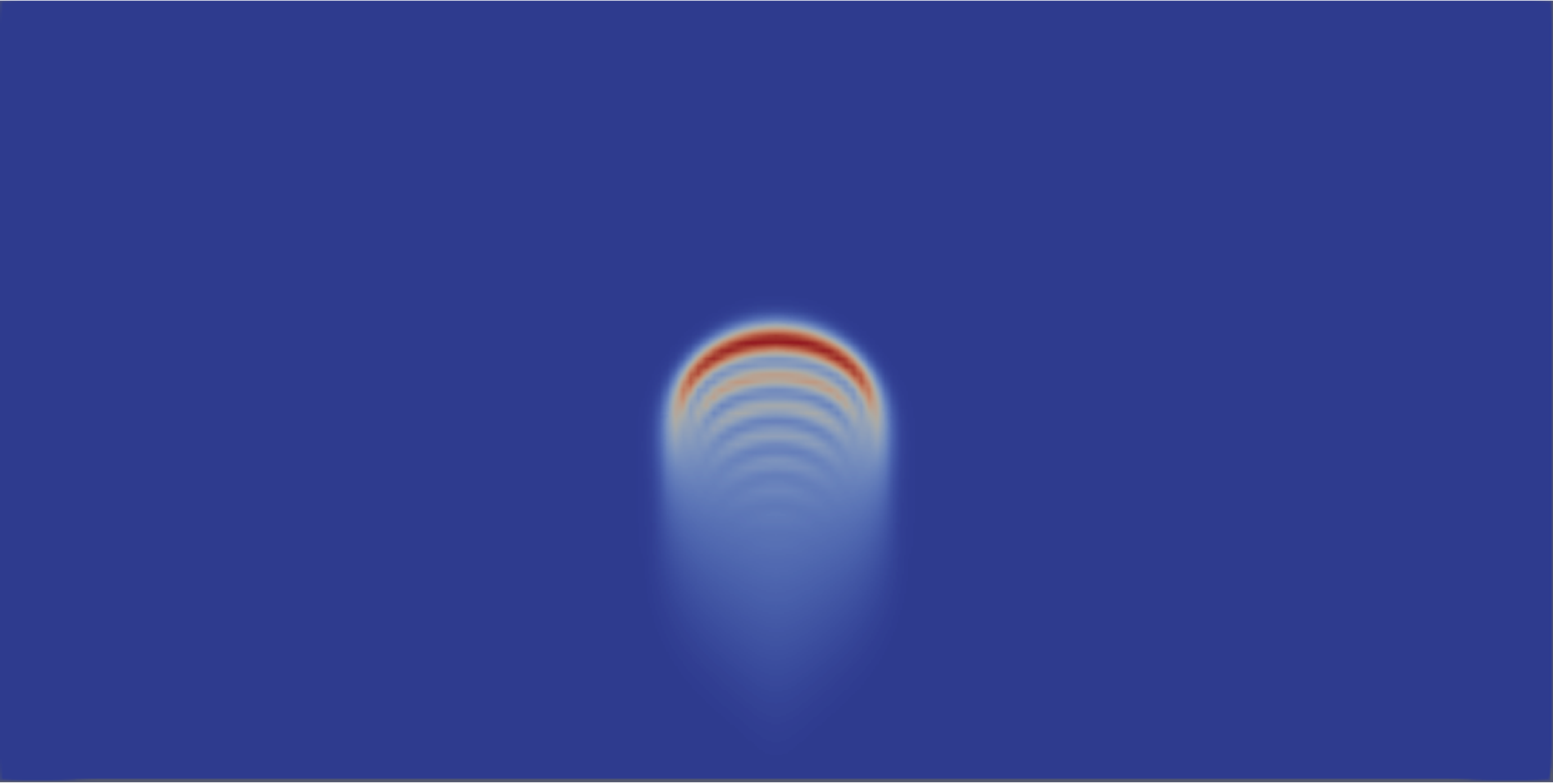}		
		\subcaption{}
		\label{fig:bub_b}
	\end{subfigure}

	\begin{subfigure}[b]{.5\linewidth}
		\centering
			\includegraphics[width=\textwidth]{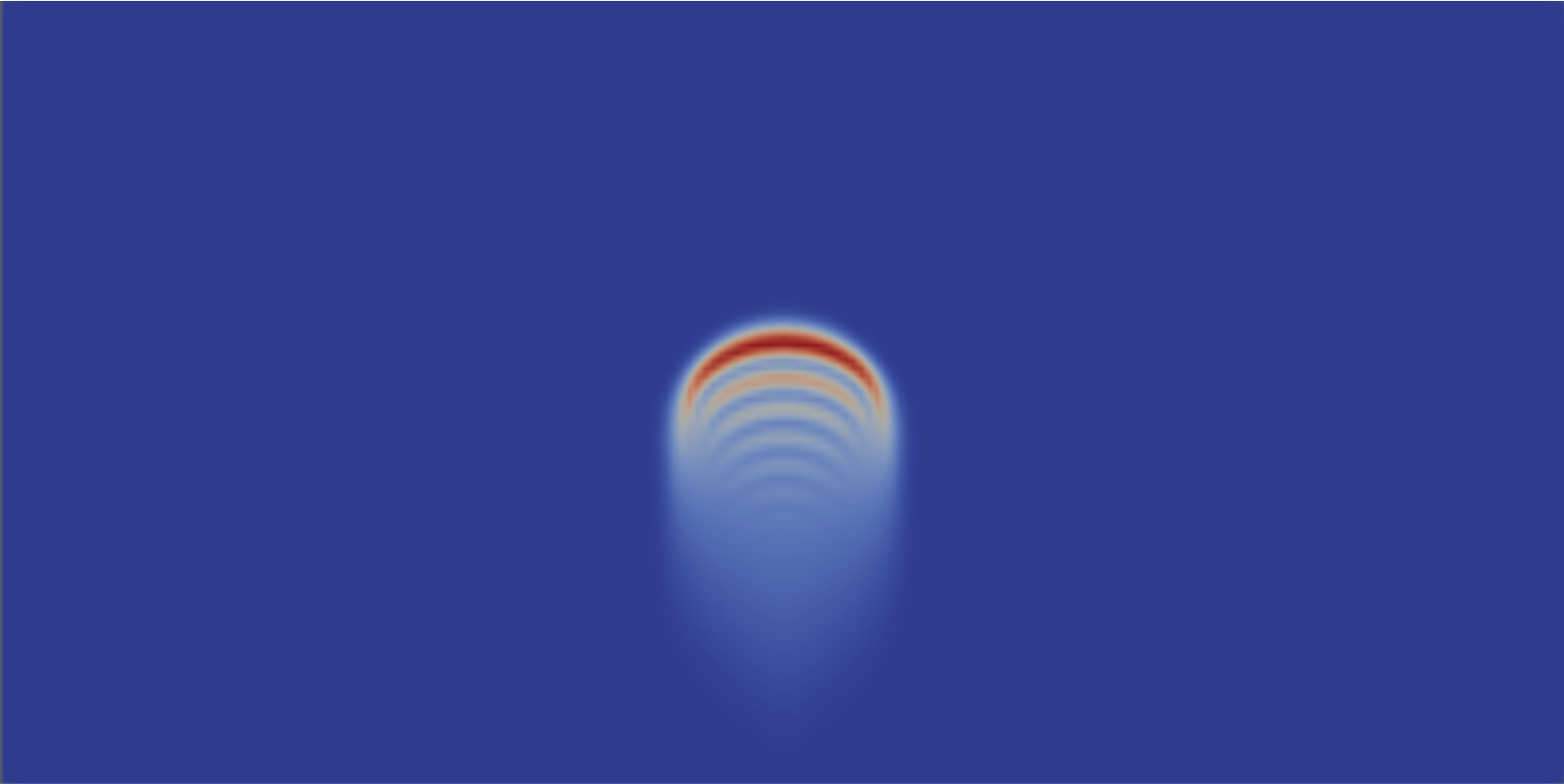}		
		\subcaption{}
		\label{fig:bub_c}
	\end{subfigure}
	\begin{subfigure}[b]{.5\linewidth}
		\centering
			\includegraphics[width=\textwidth]{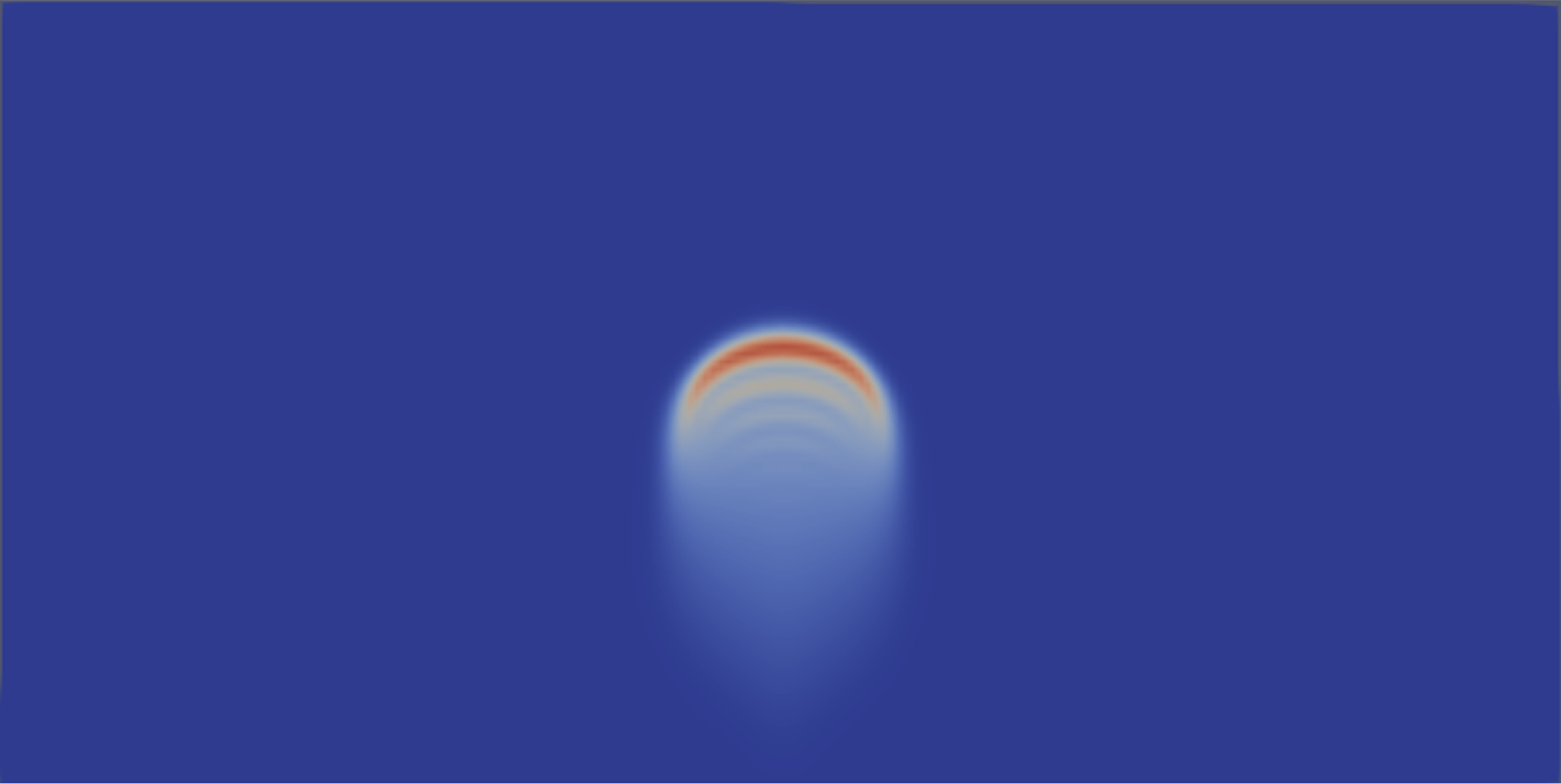}		
		\subcaption{}
		\label{fig:bub_d}
	\end{subfigure}

	\begin{subfigure}[b]{.5\linewidth}
		\centering
			\includegraphics[width=\textwidth]{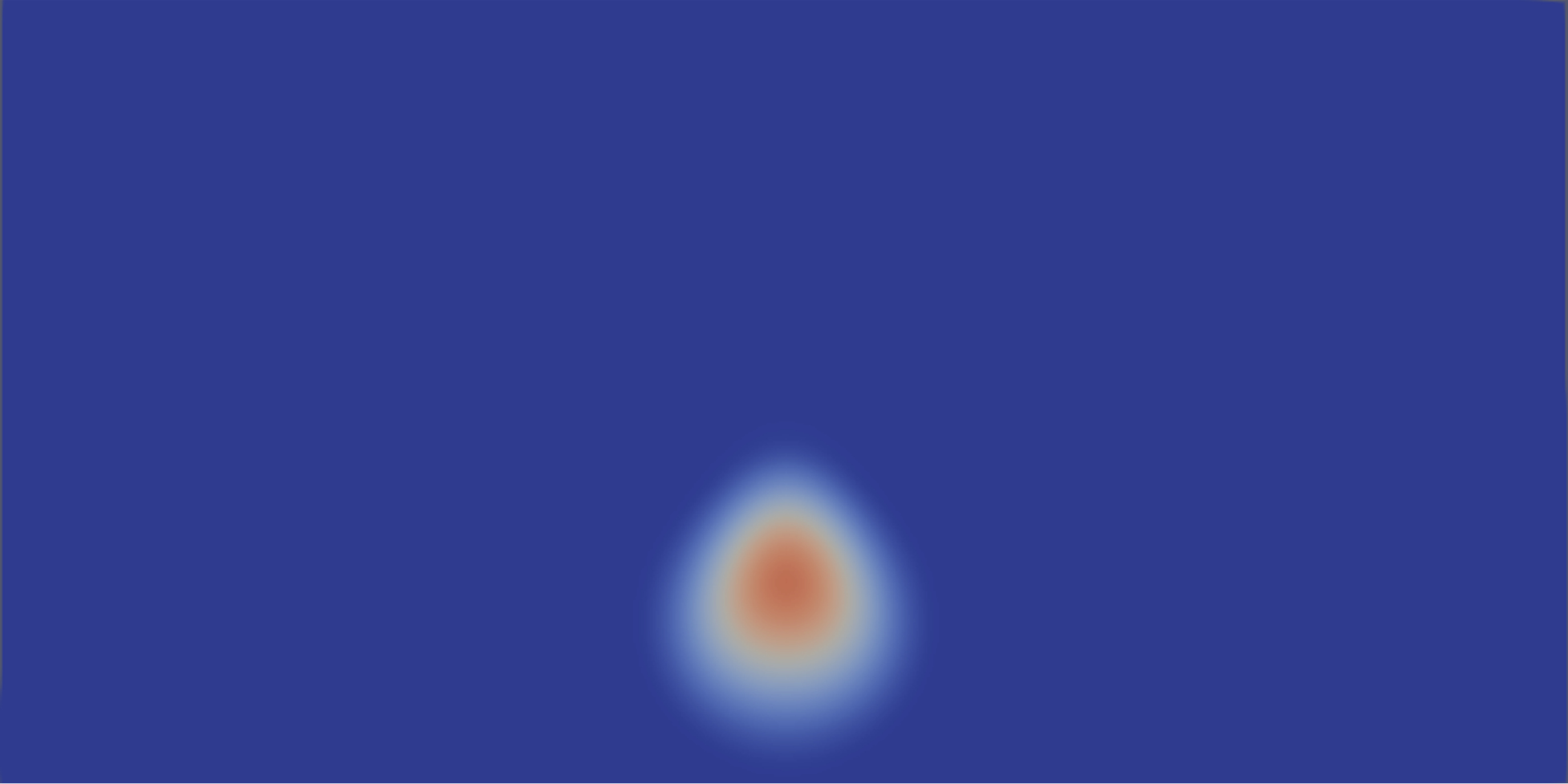}		
		\subcaption{}
		\label{fig:bub_e}
	\end{subfigure}
	\begin{subfigure}[b]{.5\linewidth}
		\centering
			\includegraphics[width=\textwidth]{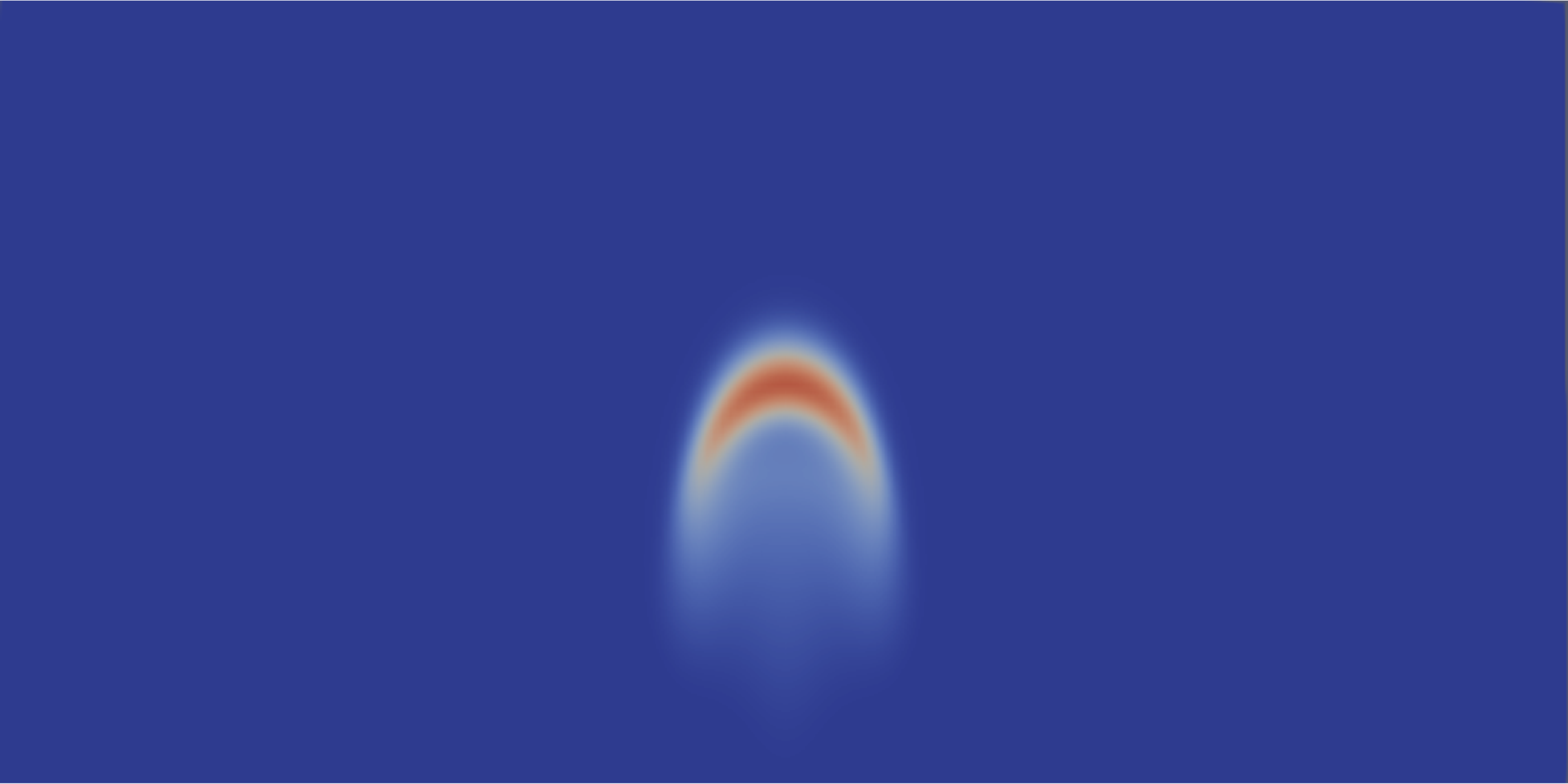}		
		\subcaption{}
		\label{fig:bub_f}
	\end{subfigure}
	
\caption{Problem {\em Thermal 1}; potential temperature distribution in a slice with a normal $z=0$ through the centre of the domain, varying between 300 K (blue colour) and 302 K (red colour) at $t=500 s$ on a grid of $100 \times 200 \times 200$ MAC cells.
(a): $\tau = 0.2 s, \mu= 1.846 \times 10^{-5} kg/m.s, K=2$; (b): $\tau = 1 s, \mu= 1.846 \times 10^{-5} kg/m.s, K=2$;
(c)  $\tau = 1 s, \mu= 1.846 \times 10^{-2} kg/m.s, K=2$; (d): $\tau = 1 s, \mu= 1.846 \times 10 kg/m.s, K=2$;
(e) $\tau = 50 s, \mu= 1.846 \times 10^{-5} kg/m.s, K=2$;  (f): $\tau = 50 s, \mu= 1.846 \times 10^{-5} kg/m.s, K=50$.
}
	\label{fig:bub1}
\end{figure}

\begin{figure}[htp] 
	\begin{subfigure}[b]{.5\linewidth}
		\centering
			\includegraphics[width=\textwidth]{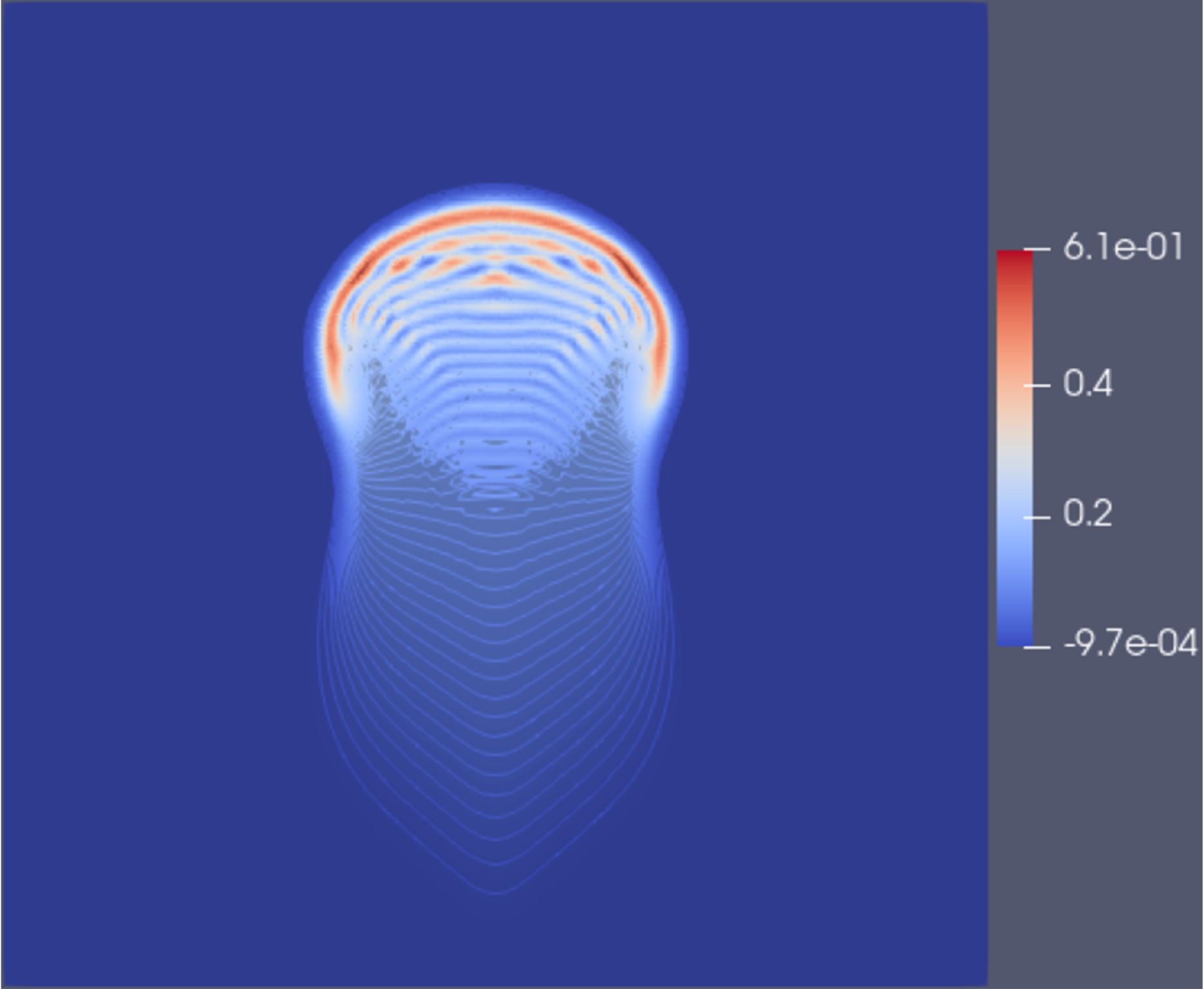}		
		\subcaption{}
		\label{fig:bub_a}
	\end{subfigure}
	\begin{subfigure}[b]{.5\linewidth}
		\centering
			\includegraphics[width=\textwidth]{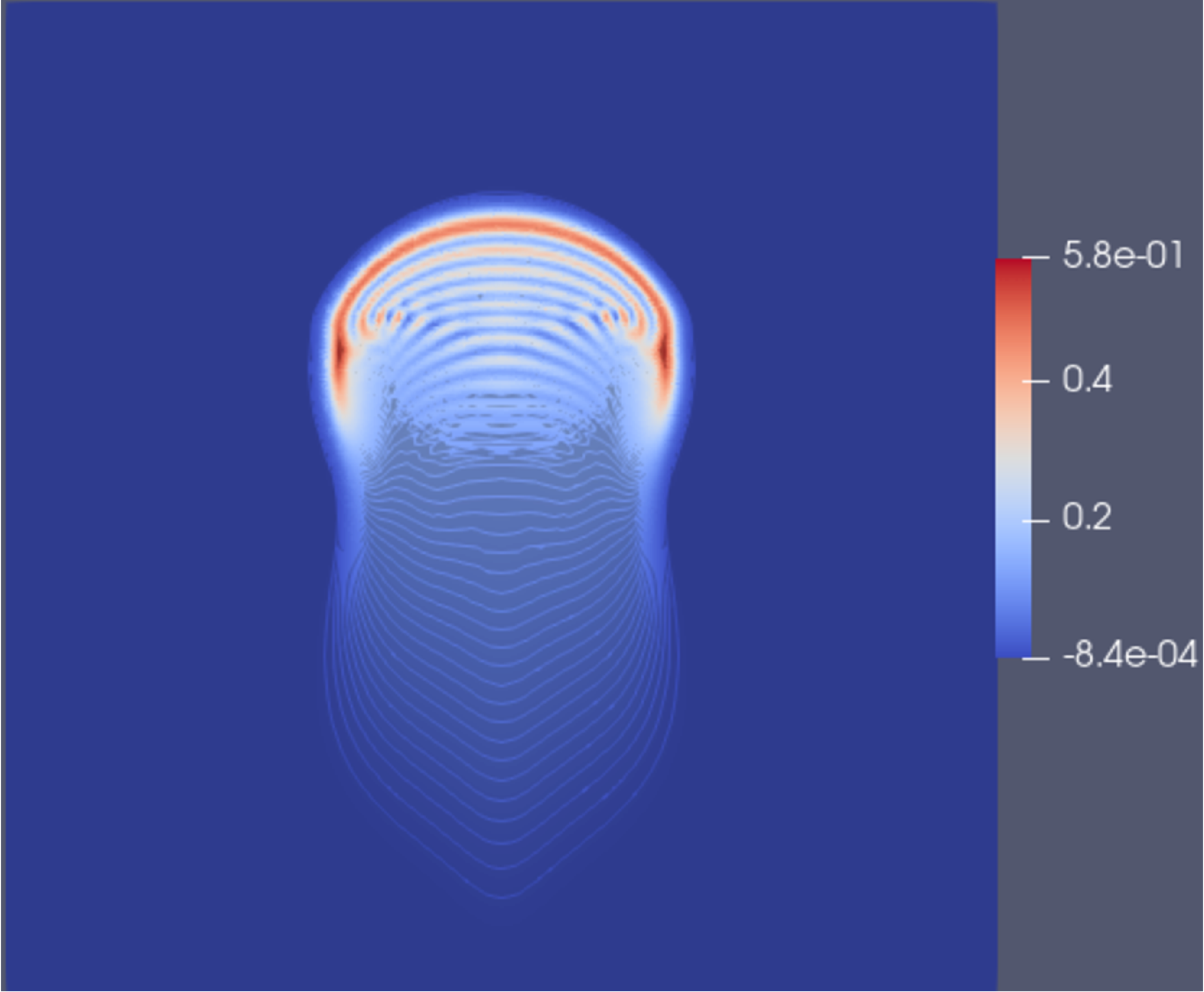}		
		\subcaption{}
		\label{fig:bub_b}
	\end{subfigure}
\begin{center}
	\begin{subfigure}[b]{.5\linewidth}
		\centering
			\includegraphics[width=\textwidth]{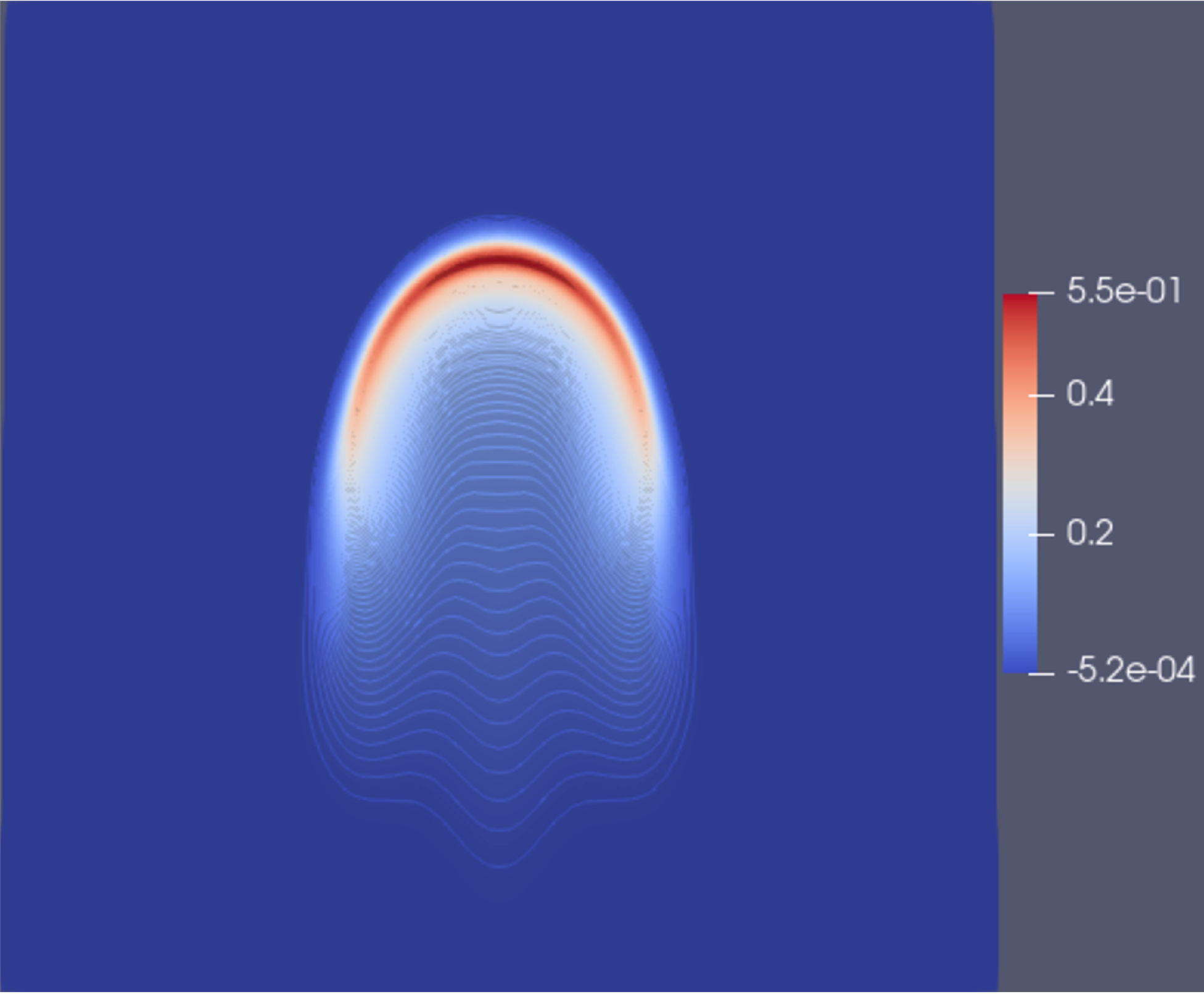}		
		\subcaption{}
		\label{fig:bub_c}
	\end{subfigure}
\end{center}
\caption{Problem {\em Thermal 2}; potential temperature perturbation in a slice with a normal $z=0$ through the centre of the domain; contour lines varying between 0 K and 0.5 K, with a step of 0.005, at $t=400 s$, on a grid of $200 \times 200 \times 200$ MAC cells.
(a): $\tau = 0.25 s, K=2$; (b): $\tau = 1 s, K=10$;
(c)  $\tau = 1 s, K=2$.
}
	\label{fig:bub2}
\end{figure}

}

\subsection{Weak scalability.}
\label{3Weak}
\begin{table}
	\begin{tabular}{ |p{2.0cm}||p{2.1cm}|p{2.2cm}|p{2.2cm}|p{2.3cm}|p{2.7cm}|  }
		\hline
		\# cores & 1(1$\times$1$\times$1) & 32(2$\times$4$\times$4) & 256(4$\times$8$\times$8) & 512(4$\times$8$\times$16)&1024(4$\times$16$\times$16)\\
		\hline
		Efficiency (one core) &  --  & $91$ \% & $85$ \%& $82$ \%& $77$  \%\\ 
		Efficiency (one node)&--& -- & $93$ \%&  $90$ \%&  $85$ \% \\
		\hline
	\end{tabular}
	\caption{Weak Scalability test, $3375 \times 10^3$ grid points per core.}
	\label{tab:scal}
\end{table}
In this section we demonstrate the parallel performance of the method by providing weak scalability results. We consider $3375 \cdot 10^3$ grid points per core and measure the efficiency on $32$ cores ($1$ computational node), $256$ cores ($8$ nodes), $512$ cores ($16$ nodes), and $1024$ cores ($32$ nodes). The efficiency results are given in Table \ref{tab:scal}. The efficiency is given relative to $1$ core and relative to $1$ node since the drop in efficiency from $1$ to $32$ cores is likely caused by the need to share the memory bandwidth and cache with a smaller number of cores within the computational node, rather than scaling properties of the method (see e.g. \cite{Keating}, p. 152). The weak scaling test demonstrates excellent parallel performance. 

The scaling tests are performed using the Compute Canada Graham cluster (see www.computecanada.ca) of $2.1$GHz Intel $E5-2683$ v4 CPU cores, $32$ cores per node, and each node connected via a $100$ Gb/s network.  The results were calculated using the wall clock time taken to simulate $10$ time steps { with two Picard iterations each. }
These computations were performed three times for each configuration, and the average wall clock time was used to compute the efficiency.  
{ On a single core, the CPU time needed to compute a problem with $3375 \times 10^3$ grid points, with one iteration, breaks down as follows: $30.5 s$ for the assembly of the right hand side, $8.1 s$ for the assembly of the matrices, $166.5 s$ for the linear solves. }

\clearpage
\section{Conclusion.}
\label{3conclusion}

The numerical experiments presented above demonstrate the effectiveness of implicit methods based on the direction splitting approach for modelling compressible flows in spherical shells in nearly incompressible and weakly compressible regimes. The proposed algorithm retains theoretically expected convergence rates and remains stable for extremely small values of the characteristic Mach number (at least as low as $M_0=10^{-6}$). The staggered spatial discretization on the MAC stencil, commonly used in numerical methods for incompressible Navier-Stokes equations, was found to be convenient for the discretization of the compressible Navier-Stokes equations written in the non-conservative form in terms of the primitive variables. This approach helped to avoid the high-frequency oscillations without any artificial stabilization terms. Nonlinear Picard iterations with the splitting error reduction were also implemented to allow one to obtain a solution of the fully nonlinear system of equations.

These results, alongside excellent parallel performance, prove the viability of the direction splitting approach in large-scale high-resolution high-performance simulations of atmospheric and oceanic flows. Possibilities for future studies and developments include research of monotonicity preserving properties of the scheme to evaluate the need for stabilization terms for flows under extreme conditions, such as high Reynolds numbers. The influence of the linearization error on the monotonicity and stability is worth investigating as well. The computational domain should be modified to represent realistic topography, and the adaptive mesh refinement is likely to be necessary for practical applications in oceanography and atmospheric sciences.

%%%%%%%%%%%%%%%%%%%%%%%%%%%%%%%%%%%%%%%%%%%%%%%%%%%%%%%%%%%%%%%%%%%%%%%%%%%%%%%%%%%%%%%%%%%%%%%%%%%%%%%%%%%%%%%
%%%%%%%%%%%%%%%%%%%%%%%%%%%%%%%%%%%%%%%%%%%%%%%%%%%%%%%%%%%%%%%%%%%%%%%%%%%%%%%%%%%%%%%%%%%%%%%%%%%%%%%%%%%%%%%
\section*{Acknowledgments}
The authors would like to acknowledge the support, under a Discovery Grant,
of the National Science and Engineering Research Council of Canada (NSERC).

This research was enabled in part by support provided by Compute Canada (www.computecanada.ca).

\appendix

\section{ Compressible Navier-Stokes equations for dry atmosphere in primitive variables.} \label{AppendixA}

The dry dynamics of Earth's atmosphere can be modeled by the compressible Navier-Stokes equations written in the conservative form in terms of density $\rho$ $[kg/m^3]$, velocity ${\bu}$ $[m/s,m/s,m/s]$ -- Cartesian or $[m/s,1/s,1/s]$ -- spherical, and the total energy per unit volume $E$ $[J/m^3]$ (see \cite{MKMMKVGHJ2016}):

\begin{equation} \label{Bconservmass}
\frac{\partial \rho}{\partial t} + \nabla \cdot (\rho {\bu}) = 0
\end{equation}

\begin{equation}\label{Bconservmom}
\frac{\partial \rho {\bu}} {\partial t} + \nabla \cdot (\rho {\bu} \otimes {\bu}) +  \nabla p + 2 \rho \left ( {\bomega} \times {\bu}  \right ) + \rho {\bg}  -  \nabla \cdot \hat{\bsigma}   = 0
\end{equation}

\begin{equation} \label{Bconservener}
\frac{\partial E}{\partial t} + \nabla \cdot ((E+p){\bu}) - \nabla \cdot \left(\frac{\mu c_p}{Pr}\nabla T + {\bu} \cdot \hat{\bsigma} \right)= 0
\end{equation}
where ${\bomega}$ $[1/s,1/s,1/s]$ is the rotational velocity of the Earth, $\hat{\bsigma}$ is the viscous stress tensor given by $$\hat{\bsigma} = \mu \left[ \left(\nabla {\bu} + (\nabla {\bu})^T \right) - \frac{2}{3} (\nabla \cdot {\bu}) \hat{I} \right],$$ ${\bg}$ $[m/s^2,m/s^2,m/s^2]$ -- Cartesian or $[m/s^2,1/s^2,1/s^2]$ -- spherical, is the sum of the true gravity and the centrifugal force, $\displaystyle c_p$ $[J/(K \cdot kg)]$, $\displaystyle c_v$ $[J/(K \cdot kg)]$, $\mu$ $[kg/(s \cdot m)]$, Pr, $\displaystyle \gamma = \frac{c_p}{c_V}$, $\pi_{\infty}$ $[Pa]$ are constant for each material. The total energy $E$ is the sum of the internal energy ($\displaystyle e = c_V T + \frac{\pi_{\infty}}{\rho}$), kinetic energy, and gravitational potential energy: 
\begin{equation}
E = \rho e + \frac{1}{2}\rho {\bu} \cdot {\bu} + \rho g r 
\end{equation}
where r $[m]$ is the radial distance from the center of the Earth. The viscous stress tensor for a Newtonian fluid is given by
\begin{equation}
\hat{\bsigma} = \mu \left[ \left(\nabla {\bu} + (\nabla {\bu})^T \right) - \frac{2}{3} (\nabla \cdot {\bu}) \hat{I} \right].
\end{equation}
Pressure is given through the Stiffened Gas Equation of State:
\begin{equation}\label{Bsgeosb}
p = (\gamma-1) \rho e - \gamma \pi_{\infty}.
\end{equation} 
The goal of this appendix is to re-write equations (\ref{Bconservmass})-(\ref{Bconservener}) in the non-conservative form in terms of the primitive variables $p$ $[Pa]$, ${\bu}$ $[m/s,m/s,m/s]$ -- Cartesian or $[m/s,1/s,1/s]$ -- spherical, and $T$ $[K]$. Then, density will be given by the following equation of state (equivalent to \ref{Bsgeosb}):
\begin{equation}
\rho = \frac{p + \pi_{\infty}}{c_V (\gamma-1)T}.
\end{equation}
Taking into account the mass conservation (\ref{Bconservmass}), the momentum conservation can be re-written in the non-conservative form as:
\begin{equation}\label{Bnonconservmomprel}
 \frac{\partial {\bu}}{\partial t} + {\bu} \cdot \nabla {\bu} + \frac{1}{\rho} \nabla p - \frac{1}{\rho} \nabla \cdot \hat{\bsigma} + {\bg} + 2 ({\bu} \times {\bomega}) = 0.
\end{equation}

In order to rewrite the energy conservation equation (\ref{Bconservener}) in a non-conservative form, we first denote:
$$ Q = - \nabla \cdot \left(\frac{\mu c_p}{Pr}\nabla T + {\bu} \cdot \hat{\bsigma} \right), $$
and note that  (\ref{Bconservener}) can be re-written as
$$ \frac{\partial E}{\partial t} + \nabla \cdot ((E+p){\bu}) + Q = 0. $$
The total energy can be expressed as:
\begin{equation*}
\begin{split}
E = \rho c_V T + \pi_{\infty} + \frac{\rho {\bu} \cdot {\bu}}{2} &+ \rho g r = 
\frac{p + \pi_{\infty}}{\gamma-1} + \pi_{\infty } + \frac{\rho {\bu} \cdot {\bu}}{2} + \rho g r =\\
\frac{p}{\gamma-1} &+ \frac{\gamma \pi_{\infty}}{\gamma-1} + \frac{\rho {\bu} \cdot {\bu}}{2} + \rho g r.
\end{split}
\end{equation*}
Substituting the last expression into (\ref{Bconservener}) gives:
\begin{equation*}
\begin{split}
\left[ \partial_t \left( \frac{p}{\gamma-1}\right) + \nabla \cdot \left(\frac{p}{\gamma-1} {\bu} \right) + \nabla \cdot (p {\bu}) \right] + 
 \left[ \partial_t \left(\frac{\pi_{\infty} \gamma}{\gamma-1} \right) + \nabla \cdot \left(\frac{\pi_{\infty} \gamma}{\gamma-1} {\bu} \right) \right] + \\
\left[ \partial_t \left( \frac{\rho {\bu} \cdot {\bu}}{2} \right) + \nabla \cdot \left(\frac{\rho {\bu} \cdot {\bu}}{2} {\bu} \right) \right] +
\left[ \partial_t (\rho g r) + \nabla \cdot (\rho g r {\bu}) \right] + Q = 0 
\end{split}
\end{equation*}
Then using (\ref{Bconservmass}) and  (\ref{Bnonconservmomprel}), and taking into account that ${\bu} \cdot {\bg} = g u_r$ and ${\bu}$ is perpendicular to ${\bomega} \times {\bu}$, after some rearrangement one obtains
\begin{equation}\label{Bforchain}
\begin{split}
 \partial_t \left( \frac{p}{\gamma-1} \right) + \nabla \cdot \left(\frac{p}{\gamma-1} {\bu} \right) + \nabla \cdot (p {\bu})  +  \partial_t \left(\frac{\pi_{\infty} \gamma}{\gamma-1} \right) +\\
\nabla \cdot \left( \frac{\pi_{\infty} \gamma}{\gamma-1} {\bu} \right) -{\bu} \cdot \nabla p + {\bu} \cdot (\nabla \cdot \hat{\bsigma}) + Q = 0
\end{split}
\end{equation}
Let $V = {\bu} \cdot (\nabla \cdot \hat{\bsigma}) + Q$. Applying the chain rule and taking into account that the  equations
 $$\partial_t \left( \frac{1}{\gamma-1} \right) + {\bu} \cdot \nabla  \left( \frac{1}{\gamma-1} \right) = 0$$ 
 and 
 $$\partial_t \left( \frac{\pi_{\infty} \gamma}{\gamma-1} \right) + {\bu} \cdot \nabla  \left( \frac{\pi_{\infty} \gamma}{\gamma-1} \right) = 0$$ 
 represent the advection of a material interface and thus have to be satisfied, (\ref{Bconservener}) can be written as:
\begin{equation}\label{Bnonconservenerprel}
\frac{\partial p}{\partial t} + {\bu} \cdot \nabla p  +\gamma (p+\pi_{\infty}) \nabla \cdot {\bu}   + (\gamma-1)V = 0.
\end{equation}

Finally, we  substitute  $$ \rho = \frac{p}{c_V (\gamma-1)T} + \frac{\pi_{\infty}}{c_V (\gamma-1)T}, $$
into (\ref{Bconservmass}) so that it becomes:
\begin{equation*}
\begin{split}
&\frac{1}{c_V (\gamma-1)T} \left[\frac{\partial p}{\partial t} + {\bu} \cdot \nabla p \right] - 
\frac{p+\pi_{\infty}}{c_V (\gamma-1)T^2} \left[\frac{\partial T}{\partial t}+ {\bu} \cdot \nabla T \right] + \frac{p+\pi_{\infty}}{c_V (\gamma-1)T} \nabla \cdot {\bu} + \\
&\frac{p}{T} \left[ \frac{\partial }{\partial t} \left(\frac{1}{c_V (\gamma-1)} \right)+ {\bu} \cdot \nabla \left( \frac{1}{c_V (\gamma-1)}\right) \right]    
    + 
\frac{1}{T} \left[ \frac{\partial }{\partial t} \left(\frac{\pi_{\infty}}{c_V (\gamma-1)} \right)+ {\bu} \cdot \nabla \left( \frac{\pi_{\infty}}{c_V (\gamma-1)}\right) \right] = 0    
\end{split}
\end{equation*}
Since again:
 $$\frac{\partial }{\partial t} \left(\frac{1}{c_V (\gamma-1)} \right)+ {\bu} \cdot \nabla \left( \frac{1}{c_V (\gamma-1)}\right) = 0$$ and $$ \frac{\partial }{\partial t} \left(\frac{\pi_{\infty}}{c_V (\gamma-1)} \right)+ {\bu} \cdot \nabla \left( \frac{\pi_{\infty}}{c_V (\gamma-1)}\right) =0$$
represent the advection of a material interface and thus have to be satisfied, and $\displaystyle \frac{\partial p}{\partial t} + {\bu} \cdot \nabla p$ can be expressed from (\ref{Bnonconservenerprel}), the mass conservation can be written in the non-conservative form as:
\begin{equation}\label{Bnonconservmassprel}
\begin{split}
\frac{\partial T}{\partial t} + {\bu} \cdot \nabla T + (\gamma-1)T \nabla \cdot {\bu} + \frac{(\gamma-1)T}{p+\pi_{\infty}}V = 0.
\end{split}
\end{equation}
Taking into account the symmetry of the stress tensor, $V$ can be expressed as:
\begin{equation}\label{Bv}
V = - \nabla \cdot \left( \frac{\mu c_p}{Pr} \nabla T \right) - \nabla {\bu} : \hat{\bsigma}.
\end{equation}
Substituting (\ref{Bv}) into (\ref{Bnonconservenerprel}) and (\ref{Bnonconservmassprel}), we finally obtain the system (\ref{Bconservmass})-(\ref{Bconservener}) in the non-conservative form in terms of the primitive variables $(p,{\bu},T)$:

\begin{equation}\label{Bnonconservmass}
\begin{split}
\frac{\partial T}{\partial t} + {\bu} \cdot \nabla T + (\gamma-1)T \nabla \cdot {\bu} - \frac{(\gamma-1)T}{p+\pi_{\infty}}&\nabla \cdot \left( \frac{\mu c_p}{Pr} \nabla T \right) - \\&\frac{(\gamma-1)T}{p+\pi_{\infty}}\nabla {\bu} : \hat{\bsigma} = 0,
\end{split}
\end{equation}
\begin{equation}\label{Bnonconservmom}
\begin{split}
\frac{\partial {\bu}}{\partial t} + {\bu} \cdot \nabla {\bu} + \frac{1}{\rho} \nabla p - \frac{1}{\rho} \nabla \cdot \hat{\bsigma} + {\bg} + 2 ({\bu} \times {\bomega}) = 0,
\end{split}
\end{equation}
\begin{equation}\label{Bnonconservener}
\begin{split}
\frac{\partial p}{\partial t} + {\bu} \cdot \nabla p  +\gamma (p+\pi_{\infty}) \nabla \cdot {\bu}    - (\gamma-1)&\nabla \cdot \left( \frac{\mu c_p}{Pr} \nabla T \right) - \\& (\gamma-1)\nabla {\bu} : \hat{\bsigma} = 0,
\end{split}
\end{equation}
where 
\begin{equation}\label{Bnonconserveos}
\begin{split}
\rho = \frac{p + \pi_{\infty}}{c_V (\gamma-1)T}.
\end{split}
\end{equation}

It is remarkable that the equation for the temperature is actually derived from the mass conservation, not the energy conservation equation.

\section{Governing equations in spherical coordinates and definition of operators.} \label{AppendixB}

Here we aim at rewriting the system (\ref{Bnonconservmass})-(\ref{Bnonconserveos}) in spherical coordinates.
Recall that the spherical transformation is given by:
\begin{align*}
\begin{cases}
 x &= r \sin \theta \cos \phi \\
 y &= r \sin \theta \sin \phi \\
 z &= r \cos \theta .
 \end{cases}
\end{align*}
In what follows, all differential operators in spherical coordinates are denoted by the corresponding symbol with a tilde above it i.e. $\tilde{\nabla}$ denotes the gradient in spherical coordinates, an so on. Let us also denote by $\pmb{e_r}$, $\pmb{e_{\theta}}$, and $\pmb{e_{\phi}}$ the unit vectors in spherical coordinates.

To simplify the notations, we denote $\displaystyle \kappa = \frac{\mu c_p}{Pr}$ and define the following operators: 
\begin{equation} \label{OpG}
\pmb{\tilde{G}}= \pmb{ \tilde{\nabla}} \bu
\end{equation}
\begin{equation} \label{OpA1}
\pmb{A_1}({\bu}) f = {\bu} \cdot \tilde{\nabla} f
\end{equation}
\begin{equation} \label{OpA2}
\pmb{A_2}(p,{\bu}) {\bv} = \gamma \left( p + \pi_{\infty} \right) \tilde{\nabla} \cdot {\bv}  - \left( \gamma - 1\right)
\tilde{\nabla} {\bu} : \hat{\bsigma}\left({\bv}\right)
\end{equation}
\begin{equation} \label{OpA3}
\pmb{A_3} f =  - \left(\gamma - 1 \right) 
\tilde{\nabla} \cdot \left(\frac{\mu c_p}{Pr} \tilde{ \nabla} f \right)
\end{equation}
\begin{equation} \label{OpB1}
\pmb{B_1}(\rho) f = \frac{1}{\rho} \tilde{\nabla} f
\end{equation}
\begin{equation} \label{OpB2}
\pmb{B_2}(\rho,{\bu}) {\bv} = {\bu} \cdot \tilde{\nabla} {\bv} + 2 \left ( {\bomega} \times {\bv}  \right ) - \frac{1}{\rho} \left( \tilde{\nabla} \cdot ( \mu \tilde{\nabla} {\bv}) +  \tilde{\nabla} \left(\frac{\mu}{3} \tilde{\nabla} \cdot {\bv} \right) \right)
\end{equation}
\begin{equation} \label{OpC2}
\pmb{C_2}(T,p,{\bu}) {\bv} = (\gamma-1) T \tilde{\nabla} \cdot {\bv} - \frac{\left(\gamma - 1 \right)T}{p + \pi_{\infty}}
\tilde{\nabla} {\bu} : \hat{\bsigma}\left({\bv} \right)
\end{equation}
\begin{equation} \label{OpC3}
\pmb{C_3}(T,p,{\bu}) f = {\bu} \cdot \tilde{\nabla} f - \frac{\left(\gamma - 1 \right)T}{p + \pi_{\infty}} 
\tilde{\nabla} \cdot \left(\frac{\mu c_p}{Pr} \tilde{\nabla} f \right)
\end{equation}
Then the system (\ref{Bnonconservmass})-(\ref{Bnonconserveos})  can be written as
\begin{align}
%\begin{equation}
 &\frac{\partial p}{\partial t}  + \pmb{A_1}({\bu}) p + \pmb{A_2}(p,{\bu}) {\bu} + \pmb{A_3} T = 0,\label{Coppressure}\\
%\end{equation}
%\begin{equation}
&\frac{\partial {\bu}} {\partial t}  + \pmb{B_1}(\rho) p + \pmb{B_2}(\rho,{\bu}) {\bu}  + {\bg} = 0,\label{Copvelocity}\\
%\end{equation}
%\begin{equation}
&\frac{\partial T}{\partial t} + \pmb{C_2}(T,p,{\bu}) {\bu} + \pmb{C_3}(T,p,{\bu}) T = 0. \label{Coptemperature}
%\end{equation}
\end{align}
The operators (\ref{OpA1}) - (\ref{OpC3}) can be split direction-wise as follows (note that the operators with the upper subindex $M$ include mixed derivatives, derivatives in staggered directions, and other terms that cannot be naturally incorporated into the direction splitting approach):
\begin{equation*}
\pmb{A_1} f = \pmb{A^r_1} f + \pmb{A^{\theta}_1} f + \pmb{A^{\phi}_1} f = u_r \frac{\partial f}{\partial r} +  \frac{u_{\theta}}{r} \frac{\partial f}{\partial \theta} + \frac{ u_{\phi}}{r \sin \theta} \frac{\partial f}{\partial \phi} 
\end{equation*}
\begin{equation*}
\pmb{A_2} {\bv} = \pmb{A^r_2} {\bv} + \pmb{A^{\theta}_2} {\bv} + \pmb{A^{\phi}_2} {\bv} + \pmb{A^M_2} {\bv}
\end{equation*}
\begin{align*}
%\begin{equation*}
%\begin{split}
\pmb{A^r_2} {\bv} =& \left(\gamma (p + \pi_{\infty}) +  \frac{2\mu (\gamma-1)}{3}  \left(  {\tilde{G}}_{r r} + {\tilde{G}}_{\theta \theta} +  {\tilde{G}}_{\phi \phi} \right) \right)  \frac{1}{r^2}\frac{\partial \left(r^2 v_r \right)}{\partial r}  - \\
&2\mu  (\gamma-1) {\tilde{G}}_{rr}  \frac{\partial v_r}{\partial r}  
-  \mu (\gamma-1) \left(2 {\tilde{G}}_{\theta \theta} + 2 {\tilde{G}}_{\phi \phi} \right) \frac{v_r}{r}\\
%\end{split}
%\end{equation*}
%\begin{equation*}
%\begin{split}
\pmb{A^{\theta}_2} {\bv} =& \left(\gamma (p + \pi_{\infty}) +  \frac{2\mu (\gamma-1)}{3}  \left(  {\tilde{G}}_{r r} + {\tilde{G}}_{\theta \theta} +  {\tilde{G}}_{\phi \phi} \right) \right)  \frac{1}{r \sin \theta}\frac{\partial \left(\sin \theta v_{\theta} \right)}{\partial \theta}  - \\
  &2 \mu (\gamma-1) \frac{{\tilde{G}}_{\theta \theta}}{r}  \frac{\partial v_{\theta}}{\partial \theta} - \mu (\gamma-1) \left(\frac{2 {\tilde{G}}_{\phi \phi}}{\tan \theta} - {\tilde{G}}_{r \theta} - {\tilde{G}}_{\theta r} \right) \frac{v_{\theta}}{r}\\
  \pmb{A^{\phi}_2} {\bv} =& \left(\gamma (p + \pi_{\infty}) +  \frac{2\mu (\gamma-1)}{3}  \left(  {\tilde{G}}_{r r} + {\tilde{G}}_{\theta \theta} -  2{\tilde{G}}_{\phi \phi} \right) \right)  \frac{1}{r \sin \theta}\frac{\partial v_{\phi} }{\partial \phi}  - \\
&\mu (\gamma-1) \left(\frac{2 {\tilde{G}}_{\phi \phi}}{\tan \theta} -\frac{{\tilde{G}}_{\phi \theta} +{\tilde{G}}_{\theta \phi} }{\tan \theta} - {\tilde{G}}_{r \phi} - {\tilde{G}}_{\phi r} \right) \frac{v_{\phi}}{r}\\
%%%
\pmb{A^M_2} {\bv} = &\left[ \pmb{A^M_{2,(r)}},\pmb{A^M_{2,(\theta)}},\pmb{A^M_{2,(\phi)}} \right]  {\bv} =-\mu (\gamma-1) \left( {\tilde{G}}_{r \theta} + {\tilde{G}}_{\theta r} \right)  \frac{\partial v_{\theta}}{\partial r}  
-  \mu  (\gamma-1) \left( {\tilde{G}}_{r \phi} +  {\tilde{G}}_{ \phi r} \right)  \frac{\partial v_{\phi}}{\partial r}   \\
&-  \mu (\gamma-1)  \left( {\tilde{G}}_{r \theta} +  {\tilde{G}}_{\theta r} \right)  \frac{1}{r} \frac{\partial v_r}{\partial \theta}    
-  \mu (\gamma-1)  \left( {\tilde{G}}_{\theta \phi} + {\tilde{G}}_{\phi \theta} \right)   \frac{1}{r} \frac{\partial v_{\phi}}{\partial \theta}  \\
 &-\mu  (\gamma-1) \left( {\tilde{G}}_{r \phi} +  {\tilde{G}}_{\phi r} \right)  \frac{1}{r \sin \theta} \frac{\partial v_r}{\partial \phi}   
-   \mu  (\gamma-1) \left( {\tilde{G}}_{\phi \theta} + {\tilde{G}}_{\theta \phi}  \right)   \frac{1}{r \sin \theta} \frac{\partial v_{\theta}}{\partial \phi} 
%\end{split}
%\end{equation*}
\end{align*}

%\begin{equation*}
\begin{align*}
\pmb{A_3} f =& \pmb{A^r_3} f + \pmb{A^{\theta}_3} f + \pmb{A^{\phi}_3} f = 
-\frac{\gamma-1}{r^2} \frac{\partial}{\partial r} \left(\kappa r^2 \frac{\partial f}{\partial r} \right) -\\& \frac{\gamma-1}{r^2 \sin \theta} \frac{\partial}{\partial \theta} \left(\kappa \sin \theta \frac{\partial f}{\partial \theta} \right) - \frac{\gamma-1}{r^2 \sin^2 \theta} \frac{\partial}{\partial \phi} \left(\kappa \frac{\partial f}{\partial \phi} \right) \\
%%%
\pmb{B_1}f =& \pmb{e_r} \left( \pmb{B_1^r}f \right) +  \pmb{e_{\theta}} \left( \pmb{B_1^{\theta}}f \right) +  \pmb{e_{\phi}} \left( \pmb{B_1^{\phi}}f \right) = \pmb{e_r} \left( \frac{1}{\rho} \frac{\partial (f)}{\partial r}  \right) +  \pmb{e_{\theta}} \left( \frac{1}{\rho r}  \frac{\partial f}{\partial \theta} \right) +  \pmb{e_{\phi}} \left(\frac{1}{\rho r \sin \theta}  \frac{\partial f}{\partial \phi}  \right),
%\end{equation*}
\end{align*}
\begin{equation*}
\begin{split}
\pmb{B_2} {\bv} &=   \pmb{e_r} \left( \pmb{B_2^{r,r}}{\bv}+\pmb{B_2^{\theta,r}}{\bv}+\pmb{B_2^{\phi,r}}{\bv} + \pmb{B_2^{M,r}}{\bv}  \right) +\\  &\pmb{e_{\theta}} \left( \pmb{B_2^{r,\theta}}{\bv} + \pmb{B_2^{\theta,\theta}}{\bv} + \pmb{B_2^{\phi,\theta}}{\bv}  + \pmb{B_2^{M,\theta}} {\bv}  \right) +\\  &\pmb{e_{\phi}} \left( \pmb{B_2^{r,\phi}}{\bv} + \pmb{B_2^{\theta,\phi}}{\bv} + \pmb{B_2^{\phi,\phi}}{\bv} + \pmb{B_2^{M,\phi}}{\bv}\right) 
\end{split}
\end{equation*}
\begin{align*}
%\begin{equation*}
%\begin{split}
\pmb{B_2^{r,r}}{\bv} =& u_r \frac{\partial v_r}{\partial r} 
- \frac{1}{\rho} \left[ \frac{1}{r^2} \frac{\partial}{\partial r} \left(\mu r^2 \frac{\partial v_r }{\partial r}\right) \right] - 
\frac{1}{3 \rho} \left(\frac{\partial}{\partial r} \left[ \frac{\mu}{r^2} \frac{\partial }{\partial r} \left(r^2 v_r \right) \right] \right)\\
%\end{split}
%\end{equation*}
\pmb{B_2^{\theta,r}}{\bv} =&  \frac{u_{\theta}}{r} \frac{\partial v_r}{\partial \theta}
- \frac{1}{\rho} \left[  \frac{1}{r^2 \sin \theta} \frac{\partial}{\partial \theta} \left(\mu \sin \theta \frac{\partial v_r}{\partial \theta} \right)   \right] \\
%%%
\pmb{B_2^{\phi,r}}{\bv} = & \frac{u_{\phi}}{r \sin \theta} \frac{\partial v_r}{\partial \phi}
- \frac{1}{\rho} \left[ \frac{1}{r^2 \sin^2 \theta} \frac{\partial}{\partial \phi } \left( \mu  \frac{\partial v_r}{\partial \phi} \right)   \right]  \\
%%%
\pmb{B_2^{M,r}}{\bv} =& \left[ \pmb{B^{M,r}_{2,(r)}},\pmb{B^{M,r}_{2,(\theta)}},\pmb{B^{M,r}_{2,(\phi)}} \right]  {\bv} = - \frac{u_{\theta} v_{\theta}}{r} - \frac{u_{\phi}v_{\phi}}{r} 
+ \frac{2 \mu v_{r}}{\rho r^2}  - 2 \omega \sin \theta v_{\phi} \\&
 - \frac{1}{3 \rho} \left( \frac{\partial}{\partial r} \left[ \frac{\mu}{r \sin \theta} \frac{\partial }{\partial \theta}\left(v_{\theta} \sin \theta \right) \right]\right) - \frac{1}{3 \rho} \left( \frac{\partial}{\partial r} \left[ \frac{\mu}{r \sin \theta} \frac{\partial v_{\phi}}{\partial \phi}  \right]\right) \\&
 - \frac{\mu}{\rho} \left[ 
- \frac{2}{r^2 \sin \theta} \frac{\partial}{\partial \theta} \left( \sin \theta v_{\theta} \right)  - \frac{2}{r^2 \sin \theta} \frac{\partial v_{\phi}}{\partial \phi} 
\right] - \\&
- \frac{1}{\rho} \left[ 
- \frac{2}{r^2 \sin \theta} \frac{\partial}{\partial \theta} \left(\mu \sin \theta v_{\theta} \right)  - \frac{2}{r^2 \sin \theta} \frac{\partial (\mu v_{\phi}) }{\partial \phi}  
\right]\\
%%%
\pmb{B_2^{r,\theta}}{\bv} =& u_r \frac{\partial v_{\theta}}{\partial r}
- \frac{1}{\rho} \left[  \frac{1}{r^2} \frac{\partial}{\partial r} \left(\mu r^2 \frac{\partial v_{\theta}}{\partial r} \right) \right] \\
%%%
\pmb{B_2^{\theta,\theta}}{\bv} =& \frac{u_{\theta}}{r} \frac{\partial v_{\theta}}{\partial \theta}
- \frac{1}{\rho} \left[ \frac{1}{r^2 \sin \theta} \frac{\partial}{\partial \theta} \left(\mu \sin \theta \frac{\partial v_{\theta}}{\partial \theta}  \right) \right] -\\& \frac{1}{3\rho} \left( \frac{1}{r^2}\frac{\partial}{\partial \theta} \left[ \frac{\mu}{ \sin \theta} \frac{\partial}{\partial \theta}   \left(v_{\theta} \sin \theta \right) \right] \right)\\
%%%
\pmb{B_2^{\phi,\theta}}{\bv} =& \frac{u_{\phi}}{r \sin \theta} \frac{\partial v_{\theta}}{\partial \phi}
- \frac{1}{\rho} \left[ \frac{1}{r^2 \sin^2 \theta} \frac{\partial}{\partial \phi }  \left( \mu \frac{\partial v_{\theta}}{\partial \phi} \right)   \right]\\
%%%
\pmb{B_2^{M,\theta}}{\bv} =& \left[ \pmb{B^{M,\theta}_{2,(r)}},\pmb{B^{M,\theta}_{2,(\theta)}},\pmb{B^{M,\theta}_{2,(\phi)}} \right]  {\bv} = \frac{u_{\theta} v_r}{r} - \frac{u_{\phi} v_{\phi}}{r \tan \theta}
+ \frac{\mu v_{\theta}}{\rho r^2 \sin^2 \theta} - 2 \omega \cos \theta v_{\phi} \\&  - \frac{1}{3 \rho} \left( \frac{1}{r^3} \frac{\partial }{\partial \theta } \left( \mu   \frac{ \partial (r^2 v_r)}{\partial r}  \right) \right) - \frac{1}{3 \rho} \left( \frac{1}{r^2}\frac{\partial}{\partial \theta} \left( \frac{\mu}{\sin \theta} \frac{\partial v_{\phi}}{\partial \phi}  \right) \right)\\&- \frac{\mu}{\rho} \left[
\frac{1}{r^2} \frac{\partial u_r }{\partial \theta}  - \frac{\cos \theta}{r^2 \sin^2 \theta} \frac{\partial u_{\phi}}{\partial \phi} 
\right] - \frac{1}{\rho} \left[-
\frac{1}{r^2} \frac{\partial (\mu u_r) }{\partial \theta}  - \frac{\cos \theta}{r^2 \sin^2 \theta} \frac{\partial (\mu u_{\phi})}{\partial \phi} 
\right] \\
%%%
\pmb{B_2^{r,\phi}}{\bv} =& u_r \frac{\partial v_{\phi}}{\partial r}
 - \frac{1}{\rho} \left[ \frac{1}{r^2} \frac{\partial}{\partial r} \left(\mu r^2 \frac{\partial v_{\phi} }{\partial r}  \right)\right] \\
 %%%
 \pmb{B_2^{\theta,\phi}}{\bu} =& \frac{u_{\theta}}{r} \frac{\partial v_{\phi}}{\partial \theta}
- \frac{1}{\rho} \left[  \frac{1}{r^2 \sin \theta} \frac{\partial}{\partial \theta} \left(\mu \sin \theta \frac{\partial v_{\phi}}{\partial \theta}  \right) \right] \\
%%%
\pmb{B_2^{\phi,\phi}}{\bv} =& \frac{u_{\phi}}{r \sin \theta} \frac{\partial v_{\phi}}{\partial \phi}
- \frac{1}{\rho} \left[\frac{1}{r^2 \sin^2 \theta} \frac{\partial}{\partial \phi } \left( \mu \frac{\partial v_{\phi}}{\partial \phi} \right) \right] - \frac{1}{3 \rho} \left( \frac{1}{r^2 \sin^2 \theta}\frac{\partial }{\partial \phi} \left( \mu \frac{\partial v_{\phi}}{\partial \phi}\right) \right)\\
%%%
\pmb{B_2^{M,\phi}}{\bv} =& \left[ \pmb{B^{M,\phi}_{2,(r)}},\pmb{B^{M,\phi}_{2,(\theta)}},\pmb{B^{M,\phi}_{2,(\phi)}} \right]  {\bv} = \frac{u_{\phi} v_{\theta}}{r \tan \theta} + \frac{u_{\phi} v_r}{r} +
 \frac{\mu v_{\phi}}{\rho r^2 \sin^2 \theta} + 2 \omega \cos \theta v_{\theta} + 2 \omega \sin \theta v_r \\&\
  - \frac{1}{3 \rho} \left( \frac{1}{r^3 \sin \theta} \frac{\partial }{\partial \phi } \left( \mu  \frac{\partial (r^2 v_r)}{\partial r}  \right) \right) - \frac{1}{3 \rho} \left( \frac{1}{r^2 \sin^2 \theta}\frac{\partial}{\partial \phi } \left( \mu  \frac{\partial ( v_{\theta} \sin \theta)}{\partial \theta} \right) \right) \\&
   - \frac{\mu}{\rho} \left[
\frac{1}{r^2 \sin \theta} \frac{\partial u_r }{\partial \phi} +  \frac{\cos \theta}{r^2 \sin^2 \theta} \frac{\partial u_{\theta}}{\partial \phi}
\right] - \frac{1}{\rho} \left[
\frac{1}{r^2 \sin \theta} \frac{\partial (\mu u_r) }{\partial \phi} +  \frac{\cos \theta}{r^2 \sin^2 \theta} \frac{\partial (\mu u_{\theta})}{\partial \phi}
\right] 
\end{align*}

\begin{equation*}
\pmb{C_2} {\bv} = \pmb{C^r_2} {\bv} + \pmb{C^{\theta}_2} {\bv} + \pmb{C^{\phi}_2} {\bv} + \pmb{C^C_2} {\bv}
\end{equation*}
\begin{align*}
%\begin{equation*}
%\begin{split}
\pmb{C^r_2} {\bv} =& \left((\gamma-1)T +  \frac{2\mu (\gamma-1)T}{3 (p + \pi_{\infty})}  \left(  {\tilde{G}}_{r r} + {\tilde{G}}_{\theta \theta} +  {\tilde{G}}_{\phi \phi} \right) \right)  \frac{1}{r^2}\frac{\partial \left(r^2 v_r \right)}{\partial r}  - \\&
2\mu  \frac{(\gamma-1)T}{p+\pi_{\infty}} {\tilde{G}}_{rr}  \frac{\partial v_r}{\partial r}  
-  \mu \frac{(\gamma-1)T}{p+\pi_{\infty}} \left(2 {\tilde{G}}_{\theta \theta} + 2 {\tilde{G}}_{\phi \phi} \right) \frac{v_r}{r}\\
%%%
\pmb{C^{\theta}_2} {\bv} =& \left((\gamma-1)T +  \frac{2\mu (\gamma-1)T}{3 (p + \pi_{\infty})}  \left(  {\tilde{G}}_{r r} + {\tilde{G}}_{\theta \theta} +  {\tilde{G}}_{\phi \phi} \right) \right)  \frac{1}{r \sin \theta}\frac{\partial \left(\sin \theta v_{\theta} \right)}{\partial \theta}   \\&
 - 2 \mu \frac{(\gamma-1)T}{p+\pi_{\infty}} \frac{{\tilde{G}}_{\theta \theta}}{r}  \frac{\partial v_{\theta}}{\partial \theta} - \mu \frac{(\gamma-1)T}{p+\pi_{\infty}} \left(\frac{2 {\tilde{G}}_{\phi \phi}}{\tan \theta} - {\tilde{G}}_{r \theta} - {\tilde{G}}_{\theta r} \right) \frac{v_{\theta}}{r}\\
%%%
\pmb{C^{\phi}_2} {\bv} =& \left((\gamma-1)T +  \frac{2\mu (\gamma-1)T}{3 (p + \pi_{\infty})} \left(  {\tilde{G}}_{r r} + {\tilde{G}}_{\theta \theta} -  2{\tilde{G}}_{\phi \phi} \right) \right)  \frac{1}{r \sin \theta}\frac{\partial v_{\phi} }{\partial \phi}  - \\&
\mu \frac{(\gamma-1)T}{p+\pi_{\infty}} \left(\frac{2 {\tilde{G}}_{\phi \phi}}{\tan \theta} -\frac{{\tilde{G}}_{\phi \theta} +{\tilde{G}}_{\theta \phi} }{\tan \theta} - {\tilde{G}}_{r \phi} - {\tilde{G}}_{\phi r} \right) \frac{v_{\phi}}{r}\\
%%%
\pmb{C^M_2} {\bv} =& \left[ \pmb{C^M_{2,(r)}},\pmb{C^M_{2,(\theta)}},\pmb{C^M_{2,(\phi)}} \right]  {\bv} =\\&
- \mu \frac{(\gamma-1)T}{p+\pi_{\infty}} \left( {\tilde{G}}_{r \theta} + {\tilde{G}}_{\theta r} \right)  \frac{\partial v_{\theta}}{\partial r}  
-  \mu  \frac{(\gamma-1)T}{p+\pi_{\infty}} \left( {\tilde{G}}_{r \phi} +  {\tilde{G}}_{ \phi r} \right)  \frac{\partial v_{\phi}}{\partial r} -  \\&
-  \mu \frac{(\gamma-1)T}{p+\pi_{\infty}}  \left( {\tilde{G}}{\tilde{G}}_{r \theta} +  {\tilde{G}}_{\theta r} \right)  \frac{1}{r} \frac{\partial v_r}{\partial \theta}    
-  \mu \frac{(\gamma-1)T}{p+\pi_{\infty}}  \left( {\tilde{G}}_{\theta \phi} + {\tilde{G}}_{\phi \theta} \right)   \frac{1}{r} \frac{\partial v_{\phi}}{\partial \theta}  \\&
-\mu  \frac{(\gamma-1)T}{p+\pi_{\infty}} \left( {\tilde{G}}_{r \phi} +  {\tilde{G}}_{\phi r} \right)  \frac{1}{r \sin \theta} \frac{\partial v_r}{\partial \phi}   
-   \mu  \frac{(\gamma-1)T}{p+\pi_{\infty}} \left( {\tilde{G}}_{\phi \theta} + {\tilde{G}}_{\theta \phi}  \right)   \frac{1}{r \sin \theta} \frac{\partial v_{\theta}}{\partial \phi} 
%\end{split}
%\end{equation*}
\end{align*}

\begin{equation*}
\pmb{C_3} f = \pmb{C^r_3} f + \pmb{C^{\theta}_3} f + \pmb{C^{\phi}_3} f 
\end{equation*}
\begin{align*}
 &\pmb{C^r_3} f =  u_r \frac{\partial f}{\partial r} 
 -\frac{(\gamma-1)T}{(p + \pi_{\infty})r^2} \frac{\partial}{\partial r} \left(\kappa r^2 \frac{\partial f}{\partial r} \right)\\
&\pmb{C^{\theta}_3} f =  \frac{u_{\theta}}{r} \frac{\partial f}{\partial \theta}
 - \frac{(\gamma-1)T}{(p + \pi_{\infty})r^2 \sin \theta} \frac{\partial}{\partial \theta} \left(\kappa \sin \theta \frac{\partial f}{\partial \theta} \right) \\
&\pmb{C^{\phi}_3} f = \frac{ u_{\phi}}{r \sin \theta} \frac{\partial f}{\partial \phi} 
- \frac{(\gamma-1)T}{(p + \pi_{\infty})r^2 \sin^2 \theta} \frac{\partial}{\partial \phi} \left(\kappa \frac{\partial f}{\partial \phi} \right).
\end{align*}
If we denote by $\bU$ the vector of unknowns: $\bU = [p,u_r,u_{\theta},u_{\phi},T]^T,$ define $\bG$ as the gravity vector $$\bG = [0,\bg^T,0]^T,$$ and combine operators in corresponding directions by introducing the following block-operators:
\begin{equation}
\label{dr}
\pmb{D_r}(\bU) = \begin{bmatrix}
\pmb{A_1^r} & \pmb{A_2^r}& 0 & 0 & \pmb{A_3^r}  \\
\pmb{B_1^r} & \pmb{B_2^{r,r}}& 0 & 0 & 0  \\
0 & 0 & \pmb{B_2^{r,\theta}}&  0 & 0  \\
0 & 0& 0 & \pmb{B_2^{r,\phi}} & 0  \\
0 & \pmb{C_2^r}  & 0 & 0 & \pmb{C_3^r}
\end{bmatrix}
\end{equation}
\begin{equation}
\label{dtheta}
\pmb{D_{\theta}} (\bU)= \begin{bmatrix}
\pmb{A_1^{\theta}} & 0& \pmb{A_2^{\theta}}  & 0 & \pmb{A_3^{\theta}}  \\
0 & \pmb{B_2^{{\theta},r}}& 0 & 0 & 0  \\
\pmb{B_1^{\theta}} & 0 & \pmb{B_2^{{\theta},\theta}}&  0 & 0  \\
0 & 0& 0 & \pmb{B_2^{{\theta},\phi}} & 0  \\
0 & 0  & \pmb{C_2^{\theta}} & 0 & \pmb{C_3^{\theta}}
\end{bmatrix}
\end{equation}
\begin{equation}
\label{dphi}
\pmb{D_{\phi}} (\bU)= \begin{bmatrix}
\pmb{A_1^{\phi}} & 0 & 0  &  \pmb{A_2^{\phi}} & \pmb{A_3^{\phi}}  \\
0 & \pmb{B_2^{{\phi},r}}& 0 & 0 & 0  \\
0 & 0 & \pmb{B_2^{{\phi},\theta}}&  0 & 0  \\
\pmb{B_1^{\phi}} & 0& 0 & \pmb{B_2^{{\phi},\phi}} & 0  \\
0 & 0  & 0 & \pmb{C_2^{\phi}} & \pmb{C_3^{\phi}}
\end{bmatrix}
\end{equation}
\begin{equation}
\label{dm}
\pmb{D_{M}} (\bU)= \begin{bmatrix}
0 & \pmb{A_{2,(r)}^M} & \pmb{A_{2,(\theta)}^M}  &  \pmb{A_{2,(\phi)}^M} & 0  \\
0 & \pmb{B_{2,(r)}^{M,r}} & \pmb{B_{2,(\theta)}^{M,r}}  & \pmb{B_{2,(\phi)}^{M,r}}  & 0  \\
0 & \pmb{B_2^{M,{2,(r)}}} & \pmb{B_{2,(\theta)}^{M,r}}  & \pmb{B_{2,(\phi)}^{M,r}}  & 0  \\
0 & \pmb{B_{2,(r)}^{M,r}} & \pmb{B_{2,(\theta)}^{M,r}}  & \pmb{B_{2,(\phi)}^{M,r}}  & 0  \\
0 & \pmb{C_{2,(r)}^M}  & \pmb{C_{2,(\theta)}^M} & \pmb{C_{2,(\phi)}^M} & 0
\end{bmatrix}
\end{equation}

the system (\ref{Coppressure})-(\ref{Coptemperature}) can be written in a compact form as:
\begin{equation}
\frac{\partial \bU}{\partial t} + \pmb{D_{r}}(\bU)\bU + \pmb{D_{\theta}}(\bU)\bU + \pmb{D_{\phi}}(\bU)\bU + \pmb{D_{M}}(\bU)\bU + \bG = 0.
\end{equation}

\section*{References}

\bibliography{refs}

\end{document}